\documentclass[cmbright,doublespace]{mmaauth}

\usepackage{moreverb}
\usepackage{graphicx} 

\usepackage[colorlinks,bookmarksopen,bookmarksnumbered,citecolor=red,urlcolor=red]{hyperref}

\newcommand\BibTeX{{\rmfamily B\kern-.05em \textsc{i\kern-.025em b}\kern-.08em
T\kern-.1667em\lower.7ex\hbox{E}\kern-.125emX}}

\newcommand{\R}{{\sf I\hspace{-.15em}R}}

\newcommand{\GP}{{\rm I\hspace{-.15em}P}}
\newcommand{\ds}{\displaystyle}
\newcommand{\qed}{\mbox{ }\nolinebreak\hfill \rule{2mm}{2mm}}

\newtheorem{theo}{Theorem}[section]
\newtheorem{rem}{Remark}[section]
\newtheorem{lem}{Lemma}[section]
\newtheorem{example}{Example}[section]

\def\volumeyear{2015}

\begin{document}

\runninghead{N.~Ovcharova}

\title{On the  coupling of regularization techniques and the boundary element method  for a hemivariational inequality 
modelling a delamination problem%\footnotemark[2]
}

\author{N.~Ovcharova %\corrauth
}

\address{ Universit\"at der Bundeswehr M\"unchen, 
         D-85577 Neubiberg/Munich, Germany }

%\corraddr{Journals Production Department, John Wiley \& Sons, Ltd,
%The Atrium, Southern Gate, Chichester, West Sussex, PO19~8SQ, UK.}

\begin{abstract}
In this paper, we couple regularization techniques of nondifferentiable optimization with the $h$-version of the boundary element method ($h$-BEM) to solve nonsmooth variational problems arising in contact mechanics. 
As a model example we consider  the delamination problem. The variational formulation of this problem leads to a hemivariational inequality (HVI) with a nonsmooth functional defined on the contact boundary. This problem is first regularized and then discretized by a $h$-BEM. We prove convergence of the $h$-BEM Galerkin solution of the regularized problem in the energy norm, provide an a-priori error estimate  and give a numerical example. 
\end{abstract}

\MOS{65M38; 74M15}


\maketitle

%\footnotetext[2]{Please ensure that you use the most up to date
%class file,
%available from the MMA Home Page at\\
%\href{http://www3.interscience.wiley.com/journal/2197/home}{\texttt{www3.interscience.wiley.com/journal/2197/home}}}

\vspace{-6pt}

\section{Introduction}
\begin{figure}[t!]
\centering
\includegraphics[trim = 0cm 0cm 0cm 0cm, clip,width=0.7\textwidth]{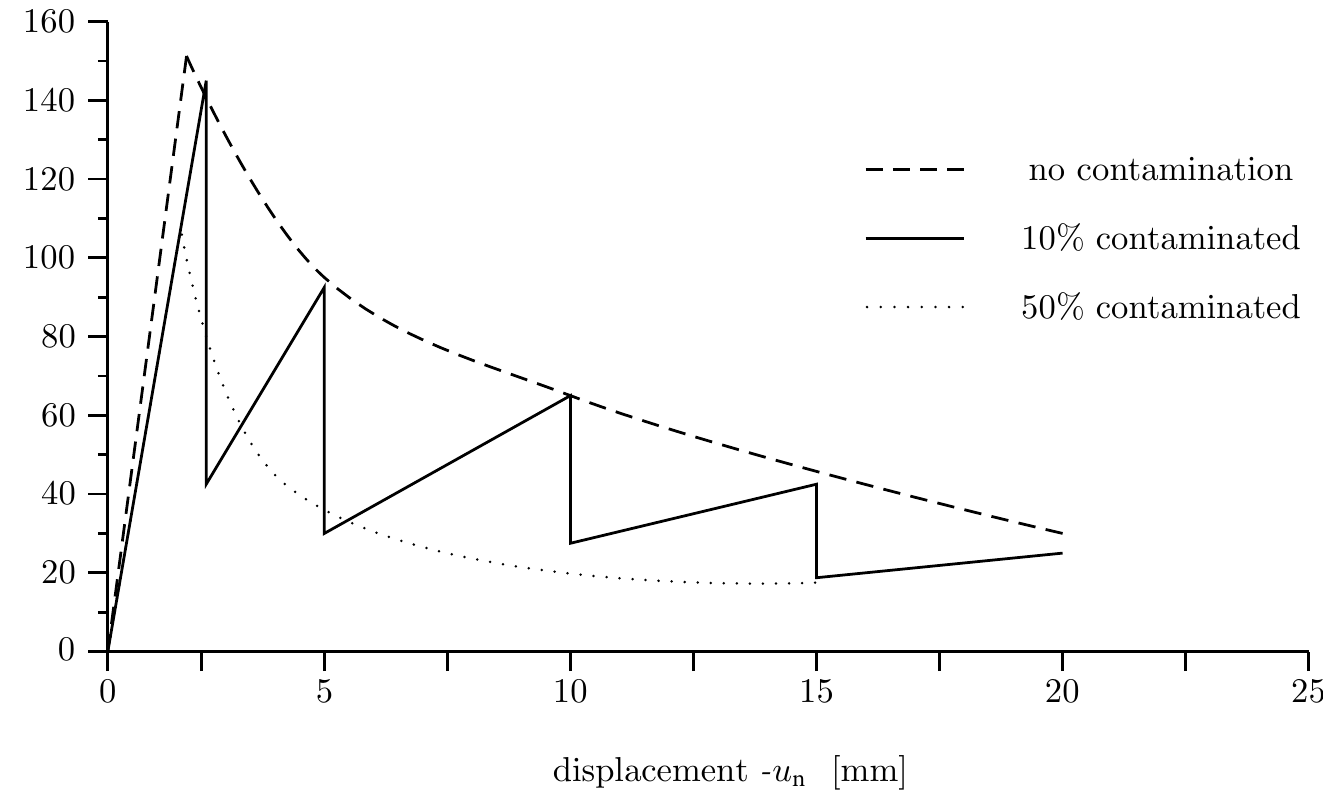}%{load-displ-curve}
\caption{Load-displacement curve determinated experimentally by DCB test for different contamination concentrations, see \cite{Gudladt}} \label{01}
\end{figure}

Efficient numerical methods for simulation of mechanical problems with nonsmooth nonmonotone contact like the adhesive contact in composite structure is of ever increasing importance in the last years. We are motivated by the delamination problems in material sciences that come from the double cantilever beam (DCB) test problem \cite{Gudladt}. The result of a typical experiment is shown in Figure \ref{01} from \cite{Gudladt}, where three probes with different levels of contamination of the interface layer have been exposed. Such problems lead in their mathematical formulation to boundary value problems involving nonmonotone and multivalued laws which can be expressed by means of the Clarke subdifferential of a nonconvex, nonsmooth  locally Lipschitz function. As a result, a nonsmooth functional defined on the contact part appears in the variational formulation of these problems. The nonsmooth behaviour in the adhesive is then modelled  by a  hemivariational inequality. There are several approaches to treat this non-differentiability. We can combine a regularization of the nonsmooth functional with finite element methods (FEM), see the PhD Thesis \cite{Ovcharova2012}, or first discretize by finite elements and then solve by nonsmooth optimization methods, %for the discrete  finite-dimensional  minimization problem, 
see \cite{ Noll_Ovcharova}. Note that in both cases, we use approximation by finite elements. Another option to treat the adhesive problem numerically is the boundary element method. To this end, the contact problem with adhesion has to be recast into a boundary integral formulation by making use of the Poincar\'{e}-Steklov operator. We emphasize that the behaviour of the adhesive interlayer is different from the contact behaviour in Signorini problems and contact problems with monotone friction.  Similar to the Coulomb friction problem, which has been treated in \cite{Eck}, the variational formulation of such problems includes a nondifferentiable functional and leads to a nonconvex problem.  For mathematical background of contact paroblems in continuum mechanics and overview of  numerical solution methods, see \cite{Gl,Gl_1981,KiOd}.  Convergence analysis and numerical solution of Signorini and friction problems by the pure $h$-boundary element Galerkin method have been discussed in \cite{Gu, Ha}. An advanced  adaptive $hp$-version of BEM for unilateral Signorini problems has been analyzed in \cite{MS-1}. In \cite{Chernov,ChSt}  new approaches based on a high-oder $hp$-BEM and a FEM-BEM coupling have been developed and applied  to provide  numerical benchmark computations for contact problems with friction. For further  numerical simulations in 2D-elasticity, we refer the reader  to \cite{cf,CoSt, MS-2}.  Multivalued boundary integral equations modelling static and dynamic contact problems have been derived and studied in \cite{Antes}. The first paper that solves HVIs modelling adhesion problems is due to Nesemann and Stephan \cite{Nesemann}. They investigate existence and uniqueness, and also propose a residual error estimator. As an exemplary function for the adhesion law they use a multivalued function with two jags. We note that their approach is based on the minimization of the potential energy function after discretization via boundary element methods with low polynomial degrees, and uses in the computations the Bundle-Newton mehod by  Luk\v{s}an and   Vl\v{c}ek   \cite{LukVl}. 

%{\color{blue}
In this paper, we focus on a contact problem with adhesive bonding and present a novel approach to solve this problem numerically, namely, we combine  regularization techniques with the $h$-BEM. More precisely, after $\varepsilon$-regularization of the nonsmooth functional, the resulting regularized problem is discretized by boundary elements. The discrete 
finite-dimensional variational inequalities  can be solved by means of numerical methods based on the optimization approach. In particular, we use an appropriate merit function to recast them into an unconstrained global minimization problem. 
We also state conditions for the uniqueness of the solution and establish $\varepsilon h$-norm convergence of the discrete solution in the energy norm. For the Galerkin solution of the regularized problem we provide an a-priori error estimate based on a novel C\'{e}a-Falk approximation lemma. The proposed approximation scheme is finally illustrated by a numerical example. Our benchmark example uses a serrated exemplary adhesion law with several jags. 
\section{ A nonmonotone boundary value problem from delamination}
\begin{figure}[b!]
\centering
\includegraphics[trim = 0cm 0cm 0cm 0cm, clip,width=0.8\textwidth]{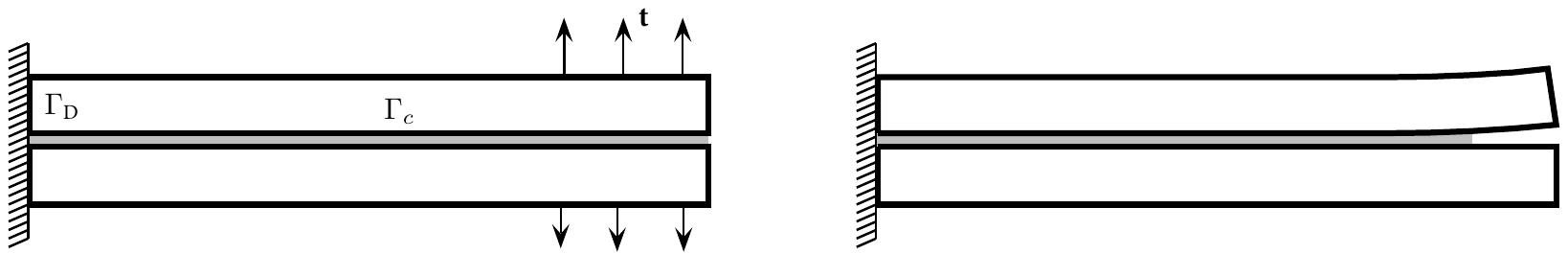}%{2Dbenchmark}
\caption{Reference congiguration for the 2D benchmark under loading. Under applied traction force $\mathbf{t}$ the crack front propagates to the left} \label{delpr}
\end{figure}

\begin{figure}[b!]
\centering
\includegraphics[trim = 5cm 15cm 0cm 0cm, clip,width=6cm]{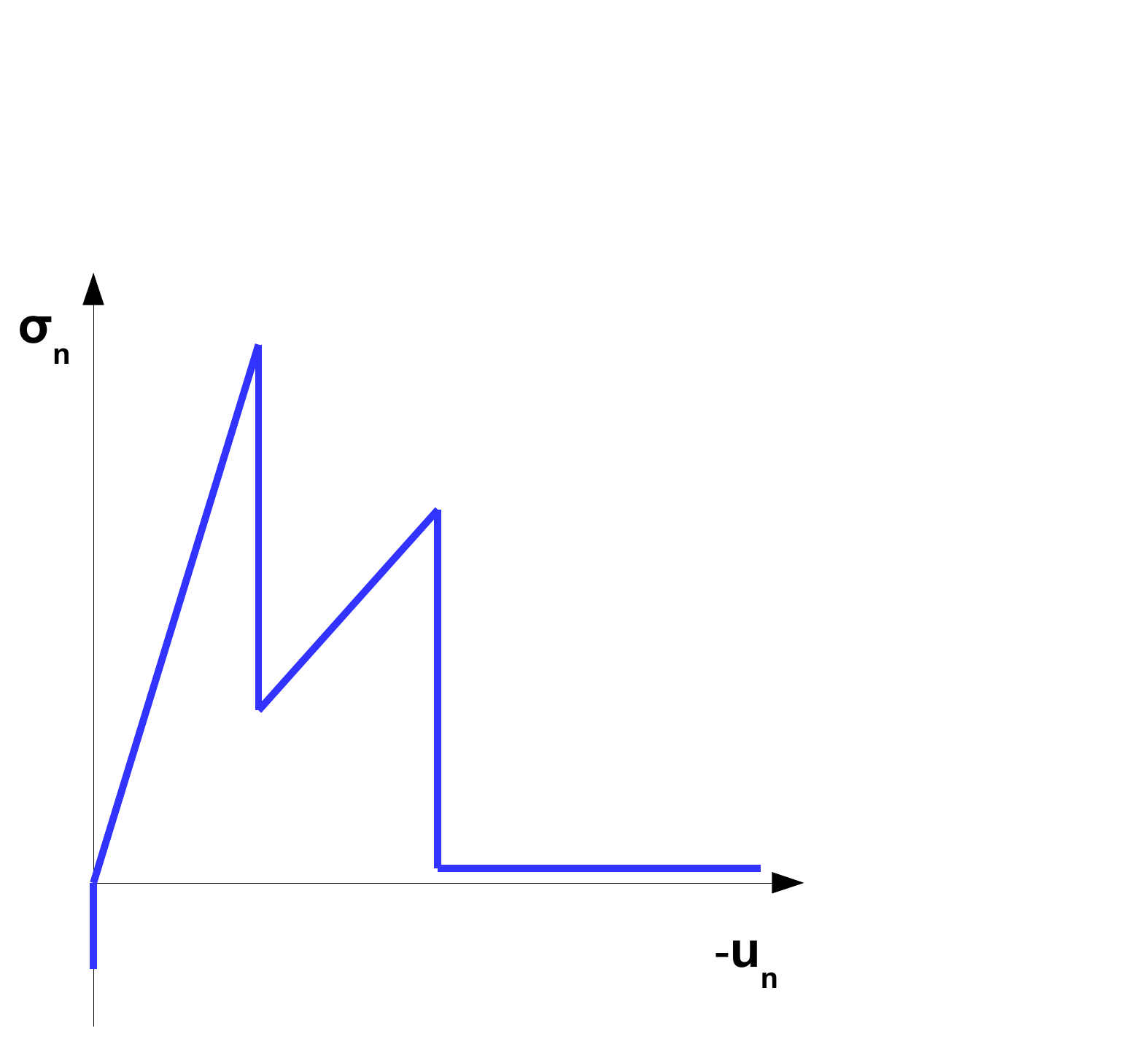}
\caption{A nonmonotone adhesion law} \label{02}
\end{figure}

Let $\Omega \subset \R^d \, (d= 2, 3)$ be a bounded  domain with Lipschitz boundary $\partial \Omega$. We assume that the boundary is decomposed into three open  disjoint parts $\Gamma_D, \Gamma_N$,
and $\Gamma_C$ such that $\partial \Omega = \overline{\Gamma}_D \cup
\overline{\Gamma}_N \cup \overline{\Gamma}_C$ and, moreover, the measures of $\Gamma_C$ and $\Gamma_D$ are positive. We consider an elastic body occupying $\overline{\Omega}$.  The body is subject to volume force ${\mathbf f} \in (L^2(\Omega))^d$. Zero displacements are prescribed on $\Gamma_D$, surface tractions $\mathbf {t}  \in (L^2(\Gamma_N))^d$ act on $\Gamma_N$, and on the part $\Gamma_C$ a nonmonotone, generally multivalued boundary condition holds.   
Further, 
$\epsilon (\mathbf {u}) = \frac{1}{2}(\nabla {\mathbf u}  + \nabla {\mathbf u}^T)$ denotes the  linearized strain tensor %(for the small deformations) to the displacement field ${\mathbf u}$ 
and $\sigma ( \mathbf {u}) = \mathcal{C}  \,: \, \epsilon (\mathbf {u})$ stands for the stress tensor, where $\mathcal{C}$ is the Hooke tensor, assumed to be 
uniformly positive definite with $L^\infty$ coefficients. The boundary stress vector  can be further decomposed into the normal, respectively, the tangential stress: 
\[
\sigma_n= \sigma(\mathbf{u}) {\mathbf n} \cdot {\mathbf n}, \quad \sigma_t=\sigma(\mathbf{u}) {\mathbf n} - \sigma_n {\mathbf n},
\]
where  ${\mathbf n}$ denotes the unit outer normal vector   on $\partial \Omega$. 
Our benchmark problem is a two- or three-dimensional symmetric laminated structure with an interlayer adhesive under loading. 
Because of the symmetry of the structure and by assuming that the forces applied to the upper and lower part of the structure are the same, it suffices to consider only the upper half of the specimen represented by $\overline{\Omega}$, see  Figure \ref{delpr} left for the  2D benchmark problem.
The delamination problem under consideration is the following. \\[0.2cm]
{\bf Problem}  $\rm{(P)}\;$  Find  ${\mathbf u} \in \mathbf{H}^1(\Omega):=[H^1(\Omega)]^d$ such that 
%\begin{equation} \label{pr1}
\begin{eqnarray}
-\mbox{div} \; \sigma ({\mathbf u})= {\mathbf f} & \; \mbox{in} \; \Omega \label{eq1} \\[0.2cm]
{\mathbf u}=0 & \;\mbox{on} \; \Gamma_D \nonumber \\[0.2cm]
\sigma({\mathbf u}){\mathbf n} = {\mathbf t} & \;\mbox{on} \; \Gamma_N \nonumber\\[0.2cm]
u_n \leq 0 & \; \mbox{on} \; \Gamma_c \nonumber\\[0.2cm]
\sigma_t({\mathbf u})=0 &  \; \mbox{on} \; \Gamma_c \nonumber\\[0.2cm]
-\sigma_n ({\mathbf u}) \in \partial f(u_n) & \; \mbox{on} \; \Gamma_c \label{nbc}
\end{eqnarray}
The contact law (\ref{nbc}), written as a differential inclusion by means  of  the Clarke subdifferential $\partial f$ \cite{Clarke} of a locally Lipschitz function $f$, describes the  nonmonotone, multivalued behaviour of the adhesive. More precisely, $\partial f$ is the physical law between the normal component $\sigma_n$ of the  boundary stress vector and the normal component $u_n= \mathbf{u} \cdot \mathbf{n}$ of the displacement $\mathbf{u}$ on $\Gamma_C$.  A typical zig-zagged nonmonotone adhesion law is shown in Figure \ref{02}.
To give a variational formulation of the above boundary value problem we define
\[
\mathbf{H}^1_D=\{\mathbf{v} \in \mathbf{H}^1(\Omega) \, : \, \mathbf{v}|_{\Gamma_D}=0 \},
\]
\[
\mathbf{K}=\{ \mathbf{v} \in \mathbf{H}^1_D \, :\, \mathbf{v}|_{\Gamma_C} \cdot \mathbf{n} \leq 0\}
\]
and introduce the bilinear form of linear elasticity
\[
a(\mathbf{u},\mathbf{v})= \int_{\Omega} \sigma (\mathbf{u}) : \varepsilon (\mathbf{v}) \, dx.
\]
Multiplying the equilibrium equation  (\ref{eq1}) in Problem $(\rm{P})$ by $\mathbf{v-u}$, integrating over $\Omega$ and applying the divergence theorem yields
%\small
\[
\int_{\Omega} \sigma (\mathbf{u}) : \varepsilon (\mathbf{v}-\mathbf{u}) \, dx = \int_{\Omega} \mathbf{f} \cdot (\mathbf{v}-\mathbf{u})\, dx + \int_{\Gamma} \sigma (\mathbf{u}) \mathbf{n} \cdot (\mathbf{v}-\mathbf{u})\, ds.
\]
From the definition of the Clarke subdifferential, the nonmonotone boundary condition (\ref{nbc}) is equivalent to 
\[
-\sigma_n (u_n)(v_n-u_n)\leq f^0(u_n;v_n-u_n) \quad \mbox{on} \; \Gamma_C.
\]
Here, the notation $f^0(x;z)$ stands for the generalized directional derivative of $f$ at $x$ in direction $z$. 

Using on $\Gamma_C$ the decomposition 
\[
\sigma (\mathbf{u})\mathbf{n} \cdot (\mathbf{v}-\mathbf{u}) = \sigma_t (\mathbf{u}) \cdot (\mathbf{v_t}-\mathbf{u_t}) +\sigma_n (\mathbf{u})(v_n-u_n)
\]
and taking into account  that on $\Gamma_C$  no tangential stresses are assumed,  
we obtain 
the  following domain hemivariational inequality: Find $\mathbf {u} \in \mathbf{K}$ such that 
\begin{equation} 
a(\mathbf{u},\mathbf{v}-\mathbf{u}) + \int_{\Gamma_C} f^0(u_n(s);v_n(s)-u_n(s)) \, ds %\nonumber \\ [0.1cm]
\geq  \int_{\Omega} \mathbf{f} \cdot (\mathbf{v}-\mathbf{u})\, dx
+ \int_{\Gamma_N} \mathbf{t} \cdot (\mathbf{v}-\mathbf{u})\, ds \quad \forall \mathbf{v} \in \mathbf{K}. \label{HemVar}
\end{equation}
We impose the following growth condition on $\partial f$: There exist positive constants $c_1$ and $c_2$ such that  for all $\xi \in \R$ and $\eta \in \partial f(\xi)$ the following inequalities hold
\begin{subequations}
 \begin{alignat}{2}
|\eta| & \leq c_1 (1+|\xi|) \label{gr1}\\
 \eta   \, \xi & \geq -c_2 |\xi| \label{gr2}
 \end{alignat}
 \end{subequations}
Note that throughout this paper $c_i$ or $C_i$ stand for  positive generic constants not necessarily the same at each occurrence.
\section{Boundary integral operator formulation}
  In this section using the Poincar\'e-Steklov operator we  rewrite the domain HVI (\ref{HemVar}) as a hemivariational inequality  defined only on the boundary. 
To this end, we introduce 
$\Gamma=\partial \, \Omega$ and $\Gamma_0=\Gamma \backslash \overline{\Gamma}_D  =\Gamma_N \cup \Gamma _C$ and define the Sobolev spaces \cite{HsWendl}:
\[
\begin{array}{lll}
H^{1/2}(\Gamma)&=&\{v\in L^2(\Gamma)\, :\, \exists v' \in H^1(\Omega), \, \mbox{tr}\, v'=v\},\\[0.2cm]
H^{1/2}(\Gamma_0)&=&\{v=v'|_{\Gamma_0} \, :\, \exists v'\in H^{1/2}(\Gamma)  \},\\[0.2cm]
\tilde{H}^{1/2}(\Gamma_0)&=&\{v=v'|_{\Gamma_0} \, : \, v'\in H^{1/2}(\Gamma)\, :\, \mbox{supp} \, v' \subset \Gamma_0 \} %\\[0.2cm]
\end{array}
\]
with the standard norms
\[
\|u\|_{H^{1/2}(\Gamma_0)}=\inf_{v\in H^{1/2}(\Gamma), v|_{\Gamma_0}=u} \|v\|_{H^{1/2}(\Gamma)}, \quad \|u\|_{\tilde{H}^{1/2}(\Gamma_0)}=\|u_0\|_{H^{1/2}(\Gamma)},
\]
where $u_0$ is the extension of $u$ by zero outside $\Gamma_0$.

The Sobolev space of negative order on $\Gamma_0$ are defined by
\[
H^{-1/2}(\Gamma_0)= (\tilde{H}^{1/2}(\Gamma_0))' \quad \mbox{and} \quad
\tilde{H}^{-1/2}(\Gamma_0) \; = \;  ({H}^{1/2}(\Gamma_0))'. 
\]
Moreover, from \cite[Lemma 4.3.1]{HsWendl} we have the inclusions
\[
\tilde{H}^{1/2}(\Gamma_0)\subset H^{1/2}(\Gamma_0)\subset L^2(\Gamma_0)\subset \tilde{H}^{-1/2}(\Gamma_0)\subset H^{-1/2}(\Gamma_0).
\]
For the spaces of vector-valued functions we use the bold symbols, e.g. 
\[
\mathbf{H}^{1/2}(\Gamma)=[H^{1/2}(\Gamma)]^d.
\]
We consider now the Navier-Lam\'{e} equation in $\R^d$, $d=2, 3$, 
\[- \mbox{div} \, \sigma(\mathbf{u})=\mathbf{f}
\]
with the Hooke's law of elasticity
\[
\sigma(\mathbf{u}) = 2 \mu \varepsilon(\mathbf{u}) + \lambda \, \mbox{div} \, \mathbf{u} \cdot \rm{I}.
\]
Here, $\rm{I}$ is the $d\times d$ identity matrix, and $\lambda, \mu >0$ are the Lam\'{e} constants depending on the material parameters: %, i.e. the modulus of elasticity $E$ and the Poisson's ratio $\nu$:
$$
\lambda = \frac{E\nu}{1-\nu^2} \,, \quad \mu= \frac{E}{1+\nu}\, . 
$$
For the solution $\mathbf{u}(\mathbf{x})$ of the Navier-Lam\'{e} equation on $\mathbf{x} \in \Omega\backslash \Gamma$ 
we have the following representation formula, also known as Somigliana's identity, see e.g. \cite{Kleiber}
\begin{equation} \label{eq2}
\hspace{-1.0cm}
\mathbf{u}(\mathbf{x})=\int_{\Gamma} E(\mathbf{x},\mathbf{y}) \left(\mathbf{T}_y \mathbf{u}(\mathbf{y})\right) ds_y -\int_\Gamma \mathbf{T}_y E(\mathbf{x},\mathbf{y})\mathbf{u}(\mathbf{y}) ds_y+ \int_\Omega E(\mathbf{x},\mathbf{y})\mathbf{f}(\mathbf{y})d\mathbf{y}, 
\end{equation}
where $E(\mathbf{x},\mathbf{y})$ is a fundamental solution of the the Navier-Lam\'{e} equation defined by
\begin{align*}
E(\mathbf{x},\mathbf{y})=\begin{cases}\displaystyle \frac{\lambda +3\mu}{4\pi \mu(\lambda+2\mu)}\left(\log |\mathbf{x}-\mathbf{y}|{\rm I} + \displaystyle \frac{\lambda + \mu}{\lambda +3\mu}\frac{(\mathbf{x}-\mathbf{y})(\mathbf{x}-\mathbf{y})^\top}{|\mathbf{x}-\mathbf{y}|^2}\right), \quad &\text{ if d=2}\\[0.2cm]
\displaystyle \frac{\lambda +3\mu}{8\pi \mu(\lambda+2\mu)}\left( |\mathbf{x}-\mathbf{y}|^{-1}{\rm I} + \displaystyle \frac{\lambda + \mu}{\lambda +3\mu}\frac{(\mathbf{x}-\mathbf{y})(\mathbf{x}-\mathbf{y})^\top}{|\mathbf{x}-\mathbf{y}|^3}\right), \quad &\text{ if d=3}
\end{cases}
\end{align*}
and $\mathbf{T}_y$  stands for the traction operator with respect to $\mathbf{y}$ defined by $\mathbf{T}_y (\mathbf{u}) := \sigma(\mathbf{u}(\mathbf{y}))\cdot \mathbf{n}_y$. 
Letting $\Omega\backslash \partial \Gamma \ni \mathbf{x} \to \Gamma$ in (\ref{eq2}), we obtain the well-known Calder\'{o}n operator
\[
\left (\begin{array}{c}
          \mathbf{u}\\[0.2cm]
         \mathbf{T}_x \mathbf{u}
         \end{array}
\right )
= \left ( \begin{array} {cc}
           \frac{1}{2}{\rm I}-K & V\\[0.2cm]
           W &      \frac{1}{2}{\rm I}+K'  
           \end{array}
   \right) 
\left( \begin{array}{c}
\mathbf{u}\\[0.2cm]
\mathbf{T}_x \mathbf{u}  
\end{array}
\right ) + \left( \begin{array}{c}
N_0\mathbf{f}\\[0.2cm]
N_1\mathbf{f}  
\end{array}
\right ),
\]   
with the single layer potential $V$, the double layer potential
$K$, its formal adjoint $K'$, and the hypersingular integral operator $W$ 
defined for $\mathbf{x} \in \Gamma$ 
as follows:
  \begin{align*}
   \left( V \phi\right ) (\mathbf{x}) & :=  \int\limits_\Gamma E(\mathbf{x},\mathbf{y}) \phi (\mathbf{y}) \, ds_y, &
   \left(K \phi \right ) (\mathbf{x}) &: =
    \int_{\Gamma} \mathbf{T}_y E^{\small T}(\mathbf{x},\mathbf{y})  \phi (\mathbf{y}) \, ds_y \\ 
   \left( K' \phi \right) (\mathbf{x}) & :=  \mathbf{T}_x \int\limits_\Gamma E 
    (\mathbf{x},\mathbf{y})  \phi (\mathbf{y}) \, ds_y, & 
    \left(W \phi \right ) (\mathbf{x}) &:= - \mathbf{T}_x  \left (K \phi\right ) (\mathbf{x}),
  \end{align*}
  and the Newton potentials $N_0, N_1$ given for $x\in \Gamma$ by
  \[
  N_0\mathbf{f}=\int_\Gamma E(\mathbf{x},\mathbf{y})\mathbf{f}(\mathbf{y}) \, ds_y, \quad N_1\mathbf{f}=\mathbf{T}_x\int_\Gamma E(\mathbf{x},\mathbf{y})\mathbf{f}(\mathbf{y}) \, ds_y. 
  \]
From  \cite{Co-1} we know that the linear operators
  \begin{align*}
    V &: \mathbf{H} ^{-1/2 + \sigma} ( \Gamma)  \rightarrow  \mathbf{H}^{1/2 + \sigma} ( \Gamma), &
     K &: \mathbf{H}^{1/2 + \sigma} ( \Gamma ) \rightarrow \mathbf{H}^{1/2 + \sigma} ( \Gamma) 
    \\
    K' &: \mathbf{H}^{-1/2 + \sigma} ( \Gamma )  \rightarrow  \mathbf{H}^{-1/2 + \sigma} (
    \Gamma), & 
    W &: \mathbf{H}^{1/2 + \sigma} ( \Gamma ) \rightarrow \mathbf{H}^{-1/2 + \sigma} (\Gamma) 
  \end{align*}
are well-defined and continuous for $|\sigma | \leq \ds\frac{1}{2}$.   
Moreover, $V$ is symmetric and positive definite (elliptic on $\mathbf{H}^{-1/2}(\Gamma)$) in $\R^3$ and, if the capacity of $\Gamma$ is smaller than 1, also in $\R^2$. This can be always  arranged by scaling, since the capacity (or conformal radius or transfinite diameter) of $\Gamma$ is smaller than 1, if $\Omega$ is 
contained in a disc with radius $<1$ 
(see e.g. \cite{SlSp,Steinbach}). The operator $W$ is symmetric and positive semidefinite with kernel $\R$ (elliptic on $\tilde{\mathbf{H}}^{1/2}(\Gamma_0)$). 
Hence, since $V$ is invertible, we obtain by taking the Schur complement  of the Calder\'{o}n projector that 
\[
\mathbf{T}_x \mathbf{u} = P\mathbf{u}-N\mathbf{f}, 
\]
where 
$P$ and $N$ are the  symmetric Poincar\'e-Steklov operator and the Newton potential given, resepectively, by
\[
P\mathbf{u} =
W \mathbf{u}  + \left( K' + \ds\frac{1}{2} I \right) V^{-1} \left( K + \ds\frac{1}{2} I \right) \mathbf{u}, \quad%=: (P\mathbf{u})(x), 
N\mathbf{f} =\left ( K'+\frac{1}{2} \mathrm{I} \right) V^{-1}N_0\mathbf{f}-N_1\mathbf{f}.
\] 
If $\mathbf{f}=0$, $P$ maps $\mathbf{u}$ to its traction and, therefore, the Poincar\'e-Steklov operator is sometimes called the Dirichlet-to-Neumann mapping.
Moreover, the operator $P$ induces a symmetric bilinear form on $\mathbf{H}^{1/2} ( \Gamma )$, and  is continuous and $\tilde{\mathbf{H}}^{1/2}(\Gamma_0)$-elliptic, 
i.e. there exist constants $c_P$, $C_P>0$ such that
\[
\begin{array}{rllll}
\|P \mathbf{u}\| _{\mathbf{H}^{-1/2}(\Gamma)}&\leq& C_P \|\mathbf{u}\| _{\mathbf{H}^{1/2}(\Gamma)} &\quad &\forall \mathbf{u} \in {\mathbf{H}^{1/2}(\Gamma)}, \\[0.2cm]
\langle P\mathbf{u}, \mathbf{u} \rangle &\geq & c_P\|\mathbf{u}\|_{\mathbf{\tilde{H}}^{1/2}(\Gamma_0)} &\quad & \forall \mathbf{u} \in {\mathbf{\tilde{H}}^{1/2}(\Gamma_0)}.
\end{array}
\]
Here, $\langle \cdot, \cdot \rangle$ is the duality pairing between the involved spaces. For the proof in 2D-case see e.g. \cite{cf}.

To simplify the notations, we introduce
$$
\mathcal{V}=\tilde{\mathbf{H}}^{1/2} (\Gamma_0), \quad  \mathcal{K}^{\Gamma}=\{ \mathbf{v} \in \mathcal{V} : \mathbf{v} |_{\Gamma_C} \cdot \mathbf{n} \leq 0 \}.
$$
Multiplying $P \mathbf{u}=\sigma(\mathbf{u})\mathbf{n}$ by  $\mathbf{v}-\mathbf{u}$, integrate on $\Gamma_0$, and using thereby again the decomposition of $\sigma (\mathbf{u})  \mathbf{n}$  on $\Gamma_C$  into the  tangential and the normal part, we get the boundary hemivariational inequality ({\bf Problem} ($\mathcal{P})$):  Find $\mathbf{u} \in \mathcal{K}^{\Gamma}$ such that 
\begin{equation} \label{BdHemVar}
\int_{\Gamma_0} (P \mathbf{u})\cdot (\mathbf{v}-\mathbf{u}) ds + 
\int_{\Gamma_C} f^0(u_n(s);v_n(s)-u_n(s)) ds %\nonumber\\ [0.1cm] 
\geq    \int_{\Gamma_0} N\mathbf{f} \cdot (\mathbf{v}-\mathbf{u}) ds + \int_{\Gamma_N} \mathbf{t} \cdot (\mathbf{v}-\mathbf{u}) ds, \quad \forall \mathbf{v} \in \mathcal{K}^{\Gamma}. %\label{BdHemVar}
\end{equation}
\normalsize
%{\color{blue}
The equivalence of  (\ref{BdHemVar}) to its corresponding domain hemivariaional problem  (\ref{HemVar}) should be understood  in the following sence. If $\mathbf{u}$ is a solution of (\ref{HemVar}), then the trace $\mathbf{u}|_\Gamma$ is a solution of the boundary hemivariational inequality (\ref{BdHemVar}). Vice versa, if $\mathbf{u}$ is a solution of (\ref{BdHemVar}) on $\Gamma$, then the extension of $\mathbf{u}$ onto $\Omega$ defined by the Somigliana's identity (\ref{eq2}) with $\mathbf{T}_x \mathbf{u}= P\mathbf{u} -N\mathbf{f}$ is a solution of (\ref{HemVar}). The same holds for the corresponding regularized problem defined in the next section. 

Finally, we note that the existence of a solution to problem (\ref{HemVar}), resp. (\ref{BdHemVar}), relies on the pseudomonotonicity of the nonsmooth boundary functional and has been investigated in \cite{Gw_Ov, Ovcharova2012, Ov_Gw}. We recall that  the functional $\varphi : X\times X \to \R$, where $X$ is a real reflexive Banach space, is pseudomonotone if $u_n\rightharpoonup u$ (weakly ) in $X$ and $\displaystyle \liminf _{n \to   \infty}  \varphi (u_n,u) \geq 0$ imply 
$\limsup _{n \to \infty} \varphi (u_n,v)\leq \varphi (u,v) $  for all $v\in X$. 
\section{Regularization of the nonsmooth functional} \label{reg}
In this section, we  recall from \cite{Ovcharova2012, Ov_Gw} a class of smoothing approximations for the maximum function based on  smoothing functions for the plus function $\mathrm{p}(t)=t^+=\max \{t, 0 \}$ and state some tools that will be used throughout this paper.

We introduce the notations 
\[
\R_{+}=\{\varepsilon\in \R \, :\, \varepsilon \geq 0\}, \quad \R_{++}=\{\varepsilon\in \R \, :\, \varepsilon > 0\}.
\]
Let 
$\hat{f} \, :\, \R_{++} \times \R$ be the smoothing function of $f$ 
defined via convolution by
\[
\hat{f}(\varepsilon, x)=\int_{\R} f(x-\varepsilon t) \rho(t) \, dt.
\] 
Here, $\varepsilon >0$ is a small regularization parameter and 
$\rho \, :\, \R \to \R_{+}$ is a probability density function such that
$$
\kappa= \int_{\R^m}|t|\, \rho(t)\, dt < \infty.
$$
In general,  the function $\hat{f}$ is not easily applicable in practice, 
but  for a special class of functions that can be expressed by means of the plus function,  
it can be explicitly computed.  For example, if $f(x)=\max\{g_1(x),g_2(x)\}$, then 
\begin{equation} \label{pr1}
f(x) = g_1(x) + \mathrm{p}[ g_2(x)- g_1(x)],
\end{equation}
  Replacing $\mathrm{p}(t)$ by its approximation $P(\varepsilon,t)$ via convulation, 
  we get  $S:\R_{++}\times \R \to \R$, 
\begin{equation} \label{smooth}
S(\varepsilon,x) =g_1(x)+ P(\varepsilon,g_2(x)-g_1(x)),
\end{equation}
as a smoothing function of $f$. 

 Using, for example, the Zang probability density function
\[
\rho(t)= \left \{ \begin{array} {ll} 1 & \mbox{if} \, - \frac{1}{2}\leq t \leq
    \frac{1}{2} \\[0.2cm]
0 & \mbox{otherwise},
\end{array} \right. 
\]
we obtain  
\begin{equation} \label{smooth_Zang}
P(\varepsilon, t) = \int _{\R} \mathrm{p}(t-\varepsilon s) \rho(s) \, ds = \left \{ \begin{array} {ll} 0 & \mbox{if} \quad t < -\frac {\varepsilon}{2}\\ [0.1cm]
\frac{1}{2\varepsilon}(t+ \frac{\varepsilon}{2})^2 & \mbox{if} \, -
\frac{\varepsilon}{2} \leq t \leq \frac{\varepsilon}{2} \\[0.1cm]
t  &\mbox{if} \quad t > \frac
    {\varepsilon}{2}
\end{array}
\right. 
\end{equation}
and hence, 
\[
S(\varepsilon, x): = 
\left \{ \begin{array} {ll} g_1(x) & \mbox{if} \, \;
    (i) \; \mbox{holds}  %\nonumber
     \\[0.1cm]
\frac{1}{2\varepsilon} [g_2(x) -g_1(x)]^2  + \frac{1}{2} (g_2(x) +
g_1(x)) + \frac{\varepsilon}{8} & \mbox{if} \, \; (ii)  \; \mbox{holds}  
\\[0.1cm] %\label{smoothZang}
g_2(x) & \mbox{if} \, \; (iii) \;  \mbox{holds}.
\end{array}
\right.
%\end{equation}
\]
The cases 
$(i)$,  $(ii)$,  $(iii)$ are defined,  respectively, by
\begin {description}
\item (i) $g_2(x) -g_1(x) \leq - \frac{\varepsilon}{2} $
\item (ii) $- \frac{\varepsilon}{2} \leq
g_2(x) -g_1(x) \leq \frac{\varepsilon}{2}$
\item (iii) $g_2(x) -g_1(x) \geq \frac{\varepsilon}{2}$.
\end{description}
For other examples of smoothing functions we refer to \cite{Ovcharova2012, Ov_Gw} and the refernces therein.  

Further, the representation formula (\ref{smooth}) can be extended to the maximum function $f:\R \to \R$ of $m$ continuous functions $g_1, \ldots, g_m$, i.e. 
\begin{align} \label{eq:bsp_f_function}
f(x)=\max \{g_1(x), g_2(x), \ldots, g_m(x)\}.
\end{align}
%of $m$ continuous functions $g_1, \ldots, g_m$. 
The smoothing function $S :\R_{++} \times \R \to R$ is then given by
\begin{equation} \label{S_general}
S(\varepsilon, x)= 
g_1(x) + P\left( \varepsilon, g_2(x)-g_1(x)+\ldots + 
P\left(\varepsilon, 
g_m(x)-g_{m-1}(x) \right)\right).
\end{equation}
The major properties of the function $S(\cdot, \cdot)$ in (\ref{S_general}) 
are listed in the following lemma:  
\begin{lem} \cite{Qi}%{Ovcharova2012}
\begin {description}
\item (i) For any $\varepsilon >0$ and for all $x\in \R$, 
\[
|S(\varepsilon,x)-f(x)|\leq (m-1)k \varepsilon.
\]
\item(ii) The function $S$ is continuously differentiable on $\R_{++}\times \R$ and for any $x\in \R$ and $\varepsilon>0$ there exist $\Lambda_i \in [0,1]$ such that $\displaystyle \sum_{i=1}^m \Lambda_i=1$ and 
 \begin {equation}\label{f1}
 \frac{\partial S(\varepsilon,x)}{\partial x} = S_x( \varepsilon, x)=\displaystyle \sum_{i=1}^m \Lambda_i  g'_i(x).
\end{equation}
Moreover,
\begin {equation}\label{f2}
\{\limsup_{z\to x, \varepsilon \to 0^+}  S_x( \varepsilon, z)\}\subseteq \partial f(x).
\end{equation}
\end {description}
\end{lem}
Assume that there exists positive constants $c_i, d_i$ such that for all $x\in \R$
\begin{subequations}
 \begin{alignat}{2}
 | g'_i(x)| &\leq c_i(1+|x|) \label{ass_1}\\ 
 g'_i(x) \, x &\geq -d_i|x|. \label{ass_2}
\end{alignat}
\end{subequations}
Under (\ref{ass_1}) - (\ref{ass_2}), the growth conditions (\ref{gr1})-(\ref{gr2}) are immediately satisfied. Moreover, from (\ref{f1})-(\ref{f2}) and (\ref{ass_1})  -(\ref{ass_2}) 
the following auxiliary result  can be easily deduced. 
\begin{lem}  
It holds that
\begin{subequations}
 \begin{alignat}{2}
  \left| S_x ( \varepsilon, x) \, z \right | &\leq  c(1+|x|)\,|z| \quad \forall x, z \in \R \label{ass_3}\\
   S_x (x, \varepsilon) \cdot (-x)  &\leq d|x| \label{ass_4}\\
  \limsup_{z\to x, \varepsilon\to 0^+} S_x(\varepsilon,x)\, \xi &\leq f^0(x;\xi) \quad \forall \xi \in \R. \label{ass_3b}
\end{alignat}
\end{subequations}
\end{lem}
Next we introduce  $J_\varepsilon :\mathbf{H}^{1/2} (\Gamma) \to \R$ defined by
$$
J_\varepsilon(\mathbf{u})=\int_{\Gamma_C}  S( \varepsilon, u_n(s)) \, ds.
$$
Since $S$ is continuously differentiable, the functional $J_\varepsilon$ is everywhere G\^{a}teaux differentiable with continuous G\^{a}teaux derivative $D J_\varepsilon : V\to V^*$ given by
\[
\langle D J_\varepsilon(\mathbf{u}), \mathbf{v} \rangle _{\Gamma_C}= \int_{\Gamma_C} S_x (\varepsilon, u_n (s)) v_n(s)  \, ds.
\]
%{\color{blue}
The {\it regularized domain problem} of (\ref{HemVar}) and the corresponding {\it regularized boundary problem}  of (\ref{BdHemVar}) are now defined, respectively,  by: 
Find $\mathbf{u}_\varepsilon \in \mathbf{K}$ such that 
   \begin{equation} \label{domainreg} 
 \hspace{-1cm}
a( \mathbf{u}_{\varepsilon}, \mathbf{v}-\mathbf{u}_{\varepsilon})  + 
\langle D J_\varepsilon(\mathbf{u}_\varepsilon), \mathbf{v} -\mathbf{u}_\varepsilon \rangle _{\Gamma_C} 
 \geq %\langle \mathbf{g}, \mathbf{v}-\mathbf{u}_\varepsilon \rangle
\int_{\Omega} \mathbf{f} \cdot (\mathbf{v}-\mathbf{u}_{\varepsilon}) \, ds+\int_{\Gamma_N} \mathbf{t} \cdot (\mathbf{v}-\mathbf{u}_{\varepsilon}) \, ds, 
\quad \forall \mathbf{v} \in \mathbf{K}, 
\end{equation}
and ({\bf Problem} ($\mathcal{P}_\varepsilon$)): Find $\mathbf{u}_\varepsilon \in \mathcal{K}^{\Gamma}$ such that 
   \begin{equation} \label{reg} 
 \hspace{-1cm}
\int_{\Gamma_0} (P \mathbf{u}_{\varepsilon})\cdot (\mathbf{v}-\mathbf{u}_{\varepsilon}) \, ds + 
\langle D J_\varepsilon(\mathbf{u}_\varepsilon), \mathbf{v} -\mathbf{u}_\varepsilon \rangle_{\Gamma_C}
 \geq %\langle \mathbf{g}, \mathbf{v}-\mathbf{u}_\varepsilon \rangle
 \int_{\Gamma_0} N\mathbf{f} \cdot (\mathbf{v}-\mathbf{u}_\varepsilon) ds + \int_{\Gamma_N} \mathbf{t} \cdot (\mathbf{v}-\mathbf{u}_{\varepsilon}) \, ds, 
 \quad \forall \mathbf{v} \in \mathcal{K}^{\Gamma}.
\end{equation}
\normalsize
According to  \cite[Theorem 4.1]{Ov_Gw} the regularized domain problem (\ref{domainreg}) has at least one solution $\mathbf{u}_\varepsilon$. Moreover,  there exists a subsequence of solutions $\{\mathbf{u}_{\varepsilon_k}\}$, $\varepsilon_k\to 0^+$,  which converges strongly in $\mathbf{H}^1(\Omega)$ to a solution of the problem (\ref{HemVar}). Because of the equivalence of the boundary variational formulations (\ref{BdHemVar}) and (\ref{reg}) to their corresponding domain variational problems (\ref{HemVar}) and (\ref{domainreg}),  we can formulate  the following result.
\begin{theo}
The regularized boundary problem (\ref{reg}) has at least one solution $\mathbf{u}_\varepsilon \in \mathcal{K}^\Gamma$. The family $\{\mathbf{u}_\varepsilon\}$ is uniformly bounded in $\mathcal{V}$. Moreover, there exists a subnet of $\{\mathbf{u}_\varepsilon\}$ which converges strongly in $\mathbf{H}^{1/2}(\Gamma)$ to a solution $\mathbf{u}$ of the boundary hemivariational inequality (\ref{BdHemVar}). 
%fueg which converges strongly to a solution u of the problem (p).
\end{theo}
\section{Uniqueness Result} \label{section4}
In this section, we give a new abstract uniqueness criteria for the solution of the  boundary hemivariational inequality. Whereas the uniqueness result of Nesemann and Stephan \cite{Nesemann} is limited to the concrete context, %this section presents a general uniqueness 
our result  exhibits the functional analytic structure. Moreover, we elaborate an example of a locally Lipschitz function that shows how the abstract uniqueness condition can be guaranteed. 

 To shorten the notations we introduce the functional $\varphi: \mathbf{H}^{1/2} (\Gamma)\times \mathbf{H}^{1/2} (\Gamma) \to \R$, 
\begin{equation} \label{varphi}
\varphi(\mathbf{u},\mathbf{v})=\displaystyle
{\int_{\Gamma_C}} f^0(u_n(s); v_n(s)- u_n(s)) \, ds 
\end{equation}
and the linear form 
\[
\langle \mathbf{g}, \mathbf{v} \rangle=\int_{\Gamma_0} N\mathbf{f} \cdot (\mathbf{v}) ds + \int _{\Gamma_N} \mathbf{t} \cdot \mathbf{v}\, ds.
\]
We assume that there exists a constant $\alpha \in [0,c_P)$ such that for any $\mathbf{u}, \mathbf{v} \in \mathcal{V}$ it holds
\begin{equation}
  \varphi(\mathbf{u},\mathbf{v})+\varphi(\mathbf{v},\mathbf{u})\leq \alpha \|\mathbf{u}-\mathbf{v}\|^2_{\mathcal{V}}. \label{unique}
\end{equation}
We have now the following abstract uniqueness result.
\begin{theo}\label{theo_0}
Under the assumption (\ref{unique}), there exists a unique solution of  problem $(\mathcal{P})$, which depends Lipschitz continuously
on $\mathbf{g} \in \mathcal{V}^*$.
\end{theo}
{\bf Proof}\, Assume that  $\mathbf{u}$, $\tilde{\mathbf{u}}$ are two solutions of $(\mathcal{P})$. Then the inequalities below hold:
$$
\langle P\mathbf{u}-\mathbf{g}, \mathbf{v} -\mathbf{u}\rangle_{\Gamma_0} +
\varphi(\mathbf{u},\mathbf{v}) \geq 0
\quad \forall \mathbf{v}\in \mathcal{K}^{\Gamma},
$$
$$
\langle P \tilde{\mathbf{u}}-\mathbf{g}, \mathbf{v} -\tilde{\mathbf{u}}\rangle_{\Gamma_0} +
\varphi(\tilde{\mathbf{u}}, \mathbf{v}) \geq 0
\quad \forall \mathbf{v}\in \mathcal{K}^{\Gamma}.
$$
Setting $\mathbf{v}=\tilde{\mathbf{u}}$ in the first inequality and
$\mathbf{v}=\mathbf{u}$ in the second one, and summing up the resulting inequalities, we get
\begin{eqnarray} \label{unique01}
\langle P\mathbf{u} -P \tilde{\mathbf{u}}, \tilde{\mathbf{u}}
-\mathbf{u} \rangle_{\Gamma_0} +
\varphi (\mathbf{u},\tilde{\mathbf{u}}) +
\varphi (\tilde{\mathbf{u}},\mathbf{u}) \geq 0.
\end{eqnarray}
We next use the coercivity of the operator $P$ and the assumption (\ref{unique}) to obtain
$$
c_P\|\mathbf{u}-\tilde{\mathbf{u}}\|_{\mathcal{V}}^2 \leq \varphi (\mathbf{u} ,\tilde{\mathbf{u}}) +
\varphi (\tilde{\mathbf{u}} ,\mathbf{u} )\leq
  \alpha  \|\mathbf{u}-\tilde{\mathbf{u}}\|^2_{\mathcal{V}}.
$$
Hence, since $\alpha \in [0, c_P)$,  if $\mathbf{u} \neq \tilde{\mathbf{u}}$ we receive a
contradiction. 

Let now $\mathbf{g}_i \in \mathcal{V}^*$ and denote $\mathbf{u}^i=\mathbf{u}_{\mathbf{g}_i}, \; i=1,2.$  
Analogously to (\ref{unique01}), we find that
$$
\langle P\mathbf{u}^1 -\mathbf{g}_1-P \mathbf{u}^2+\mathbf{g}_2, \mathbf{u}^2
-\mathbf{u}^1  \rangle_{\Gamma_0}  +
\varphi (\mathbf{u}^1 ,\mathbf{u}^2) +
\varphi (\mathbf{u}^2 ,\mathbf{u}^1 ) \geq 0.
$$
Hence,
$$
c_P\|\mathbf{u}^1 -\mathbf{u}^2 \|^2_{\mathcal{V}}\leq \varphi (\mathbf{u}^1 ,\mathbf{u}^2 ) +
\varphi (\mathbf{u}^2 ,\mathbf{u}^1) + \langle \mathbf{g}_1-\mathbf{g}_2, \mathbf{u}^2-\mathbf{u}^1 \rangle
$$
and by (\ref{unique}), 
$$
(c_P-\alpha)\|\mathbf{u}^1-\mathbf{u}^2\|^2_{\mathcal{V}}\leq  \langle \mathbf{g}_1-\mathbf{g}_2, \mathbf{u}^2-\mathbf{u}^1 \rangle\leq \|\mathbf{g}_1-\mathbf{g}_2\|_{\mathcal{V}^*} \|\mathbf{u}^1 -\mathbf{u}^2 \|_{\mathcal{V}}. 
$$
Also, since $\alpha <c_P$ we deduce that
$$
\|\mathbf{u}^1 -\mathbf{u}^2 \|_{\mathcal{V}} \leq \frac{1}{c_P-\alpha} \, \|\mathbf{g}_1-\mathbf{g}_2\|_{\mathcal{V}^*},
$$
which concludes the proof of the theorem.
\qed

Further, we present  a class of locally Lipschitz functions for which (\ref{unique}) is satisfied. We assume the following  so-called  {\it one-sided Lipschitz condition} on $\partial f$. Let $f:\R \to \R$ be  a function such that
\begin{equation} \label{clmont}
( \xi^*-\eta^*)\,  (\xi- \eta)  \geq -\alpha |\xi -\eta |^2 \quad 
\forall \xi^*\in \partial f(\xi), \; \forall \eta^*\in \partial f(\eta)
\end{equation}
for any $\xi, \eta \in \R$ and  some  $\alpha \geq 0$.
From the definition of the Clarke generalized derivative \cite{Clarke} we get
$$
f^0(\xi;\eta-\xi)= \max _{\xi^*\in \partial f(\xi)}  \xi^* \, (\eta -\xi). 
$$
Rewriting (\ref{clmont}) as
$$
\xi^* \, (\eta- \xi)  +  \eta^* \,  (\xi- \eta) \leq \alpha |\xi -\eta |^2
$$
we find 
$$
f^0(\xi;\eta-\xi)+ f^0(\eta;\xi-\eta) \leq \alpha |\xi-\eta|^2.
$$
Hence, using also the continuity properties of the mapping $\mathbf{u} \cdot \mathbf{n} : \mathbf{H}^{1/2}(\Gamma) \to L^2(\Gamma_C)$, we obtain
\begin{eqnarray*}
 \varphi(\mathbf{u},\mathbf{v})+\varphi(\mathbf{v},\mathbf{u})&& \hspace{-0.3cm}= 
\int_{\Gamma_C} f^0( u_n;  v_n - u_n )\, ds + \int_{\Gamma_C} f^0( v_n;  u_n - v_n )\, ds \nonumber\\[0.2cm] \qquad && 
\leq \alpha \| u_n -  v_n\|^2_{L^2(\Gamma_C)}\leq \alpha \|\mathbf{u}-\mathbf{v}\|^2_{\mathbf{H}^{1/2}(\Gamma)}.
\end{eqnarray*}
Hence, (\ref{unique}) is satisfied
provided that $\alpha \geq 0$ is sufficiently small ($\alpha < c_P$). 
\begin{rem} If $S(\varepsilon, \cdot) : \R \to \R$ satisfies (\ref{clmont}), i.e. there exists a constant $\alpha \geq 0$ such that
\begin{equation}\label{uniq_reg}
(S_x(\varepsilon, x_1)-S_x(\varepsilon, x_2))(x_1-x_2)\geq -\alpha |x_1-x_2|^2 \quad \forall x_1, \,x_2 \in \R,
\end{equation}
then the regularized problem ($\mathcal{P}_\varepsilon$) is unique solvable provided that $\alpha < c_p$. 
\end{rem}
We finish this section with a simple example for a locally Lipschitz function $f$, for which (\ref{clmont}) holds. 
\begin{example} 
Let the graph $g$ of  $\partial f$ consists of several decreasing straight line segments with negative slopes $-\alpha_i$, $i=1,\ldots,\mathcal{I}$, and nonnegative jumps  (see Figure \ref{positivejumps}), i.e. 
$$
\partial f(x) =\left[\underline{g}(x), \overline{g}(x) \right],
$$
where
$$
\underline{g}(x):=f(x-0)= \lim _{h\to 0^-}\frac{f(x+h)-f(x)}{h},
$$
$$
\overline{g}(x):= f(x+0)=\lim _{h\to 0^+}\frac{f(x+h)-f(x)}{h}
$$
and $\underline{g}(x)\leq \overline{g}(x)$. 

\begin{figure}[b] \label{positivejumps}
\centering
\includegraphics[trim=0cm 17cm 0cm 5cm,clip,width=0.8\textwidth]{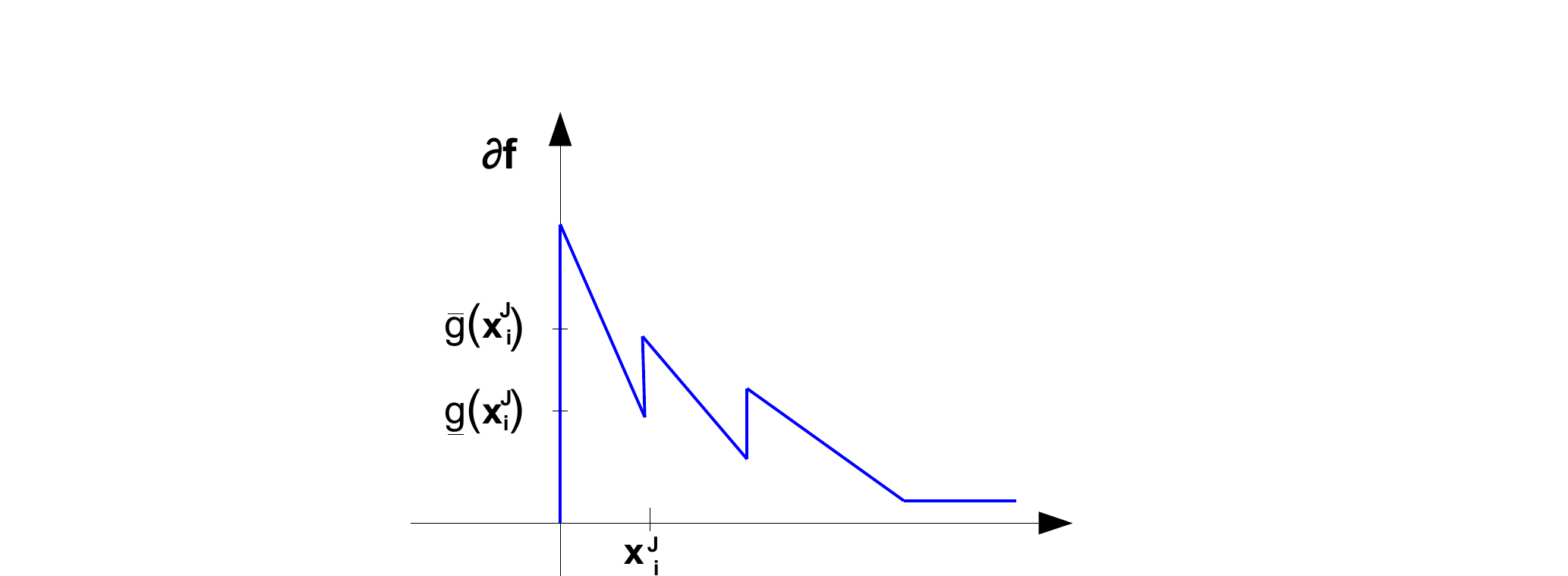} 
\caption{An example of graph of $\partial f$ with nonnegative jumps}
 \label{positivejumps}
\end{figure}
Let $-\alpha<0$ be  the  slope of the steepest decreasing segment of $\partial f$, i.e. $-\alpha=\min \{-\alpha_i : i=1,\ldots, \mathcal{I}\}$. 
Then, 
for any $x_1>x_2$ we have
\begin{eqnarray*}
\frac{\underline{g}(x_1)-\overline{g}(x_2)}{x_1-x_2} &= & \frac{\underline{g}(x_1)-\overline{g}(x_2)+ \displaystyle \sum _{i=1}^k\overline{g}(x^J_i)-\displaystyle \sum _{i=1}^k\overline{g}(x^J_i)}{x_1-x_2} \\
&\geq& \frac{\underline{g}(x_1)-\overline{g}(x_2)+ \displaystyle \sum _{i=1}^k\underline{g}(x^J_i)-\displaystyle \sum _{i=1}^k\overline{g}(x^J_i)}{x_1-x_2}\\
&= &\frac{\underline{g}(x_1)-\overline{g}(x^J_k) -\overline{g}(x_2)+\underline{g}(x^J_1)
+ \displaystyle \sum _{i=2}^k \big(\underline{g}(x^J_i)-\overline{g}(x^J_{i-1})\big)}{x_1-x_2} \\
&\geq& \frac{-\alpha(x_1-x^J_k)-\alpha(x^J_1-x_2)-\alpha \displaystyle \sum _{i=2}^k  (x^J_i -x^J_{i-1})}{x_1-x_2}\\ &=&
\frac{-\alpha(x_1-x_2)}{x_1-x_2}=-\alpha,
\end{eqnarray*}
 from which  the assumption (\ref{clmont}) follows immediately. Here, 
 $\{ x_i^J \}_{i=1}^k$ is the set of jags 
 % $k$ stands for  the number of the jags 
  between $x_1$ and $x_2$, i.e. $ x_2 < x_1^J < \cdots < x_k^J < x_1$, where $\overline{g}(x_i^j)>\underline{g}(x_i^j)$.
  %the set $\{x_i^J\}$ is ordered by $i=1,\ldots,k$. }
\end{example}

\section{Discretization with boundary elements}
Let $\Omega\subset \R^{d}$, $d=2,3$,  be a bounded domain with Lipschitz boundary $\Gamma$. As already mentioned, we only need a mesh  on the boundary.  The elements of this mesh are edges in 2D-case and triangles in 3D-case. 
 
For the discretization of the displacement $\mathbf{u}$ we use continuous piecewise linear functions on a triangulation $\mathcal{T}_h$ on $\Gamma$, which is consistent with the decomposition of $\Gamma$ into $\Gamma_0$ and $\Gamma_D$ and define
\[
\mathcal{V}_h=\{\mathbf{v_h} \in \mathbf{C}(\Gamma)\, :\, \mathbf{v_h}|_{E} \in [\GP_1 (E)]^{d-1} \quad \forall E \in \mathcal{T}_h, \;  \mathbf{v_h} =0 \; \mbox{on} \;  \overline{\Gamma}_D \} \subset \mathbf{H}_D^{1/2}(\Gamma),
\]
\[
\mathcal{K}_h^{\Gamma} =\{ \mathbf{v_h} \in \mathcal{V}_h \, :\, (\mathbf{v_h}\cdot \mathbf{n})(P_i) \leq 0 \quad \forall P_i \in \Sigma_h, \, P_i\in \overline{\Gamma}_C \backslash  \overline{\Gamma}_D\},
\]
where $\Sigma_h$ is the set of all nodes of $\mathcal{T}_h$. 

To approximate the stresses we take as ansatz space the space of piecewise constant functions on $\mathcal{T}_h$: 
\[
\mathcal{W}_h= \{ \psi \in \mathbf{L}^2(\Gamma) \, : \, \psi|_{E}  \in [\GP_0 (E)]^{d-1} \quad \forall E \in \mathcal{T}_h\} \subset \mathbf{H}^{-1/2}(\Gamma).
\]
For more details on the approximation techniques based on boundary element method see e.g. \cite{ccjg, Co-1,CoSt,Gu,GS-1,Ha,MS-1,MS-2,SlSp}. 

Let $\{\varphi_i\}_{i=1}^{N_D}$ and $\{\psi_j\}_{j=1}^{N_N}$ be the bases in $\mathcal{V}_h$ and $\mathcal{W}_h$, respectively. Then the  boundary matrices are given by
\begin{align*}
(V_{h})_{i,j}  &= \langle V \psi_i, \psi_j\rangle, \quad  (K_{h})_{i,j} = \langle K \phi_i, \psi_j\rangle, \\
(W_{h})_{i,j} &= \langle W \phi_i, \phi_j\rangle, \quad
(I_{hp})_{i,j}  = \langle \phi_i, \psi_j\rangle 
\end{align*}
The matrix $V_h$ is symmetric and positive definite, so it can be inverted by a Cholesky decomposition and as a approximation of the Galerkin matrix we obtain the matrix
\[
{P}_h=W_h+\left( K_h+ \frac{1}{2}I_h\right)^\top V^{-1}_h \left( K_h+ \frac{1}{2}I_h\right).
\]
With the canonical embeddings
\[
\begin{array}{rrlrl}
i_{h}& : & \mathcal{W}_h & \hookrightarrow& \mathbf{H}^{-1/2}(\Gamma)\\
j_{h}& : & \mathcal{V}_h & \hookrightarrow& \mathbf{H}^{1/2}(\Gamma)\\
\end{array}
\]
and their duals $i_h^*$ and $j_h^*$, 
the discrete Poincar\'{e}-Steklov operator $P_h : \mathcal{V}_h \to \mathcal{V}_h^*$ can be  also represented   by 
\[
P_h= j_h^* W j_h+ j_h^* \left(K'+\frac{1}{2}{\rm I}\right ) i_h (i_h^* V i_h)^{-1}  i_h^* \left(K+\frac{1}{2} {\rm I}\right) j_h. 
\]
According to \cite{cc}, there exists a constant $c>0$ such that 
\begin{equation} \label{coerc}
\langle P_h \mathbf{u}_h, \mathbf{u}_h \rangle_{\Gamma_0} \geq c \|j_h \mathbf{u}_h\|_{\mathbf{\tilde{H}}^{1/2}(\Gamma_0)} \quad \forall \mathbf{u}_h \in \mathcal{V}_h.
\end{equation}
Further, we define  the operator 
$E_h: \mathbf{H}^{1/2}(\Gamma)\to {\mathbf{H}}^{-1/2}(\Gamma)$, reflecting the consistency error in the discretization of the Poincar\'{e}-Steklov operator $P$, by 
\[
E_h:=P-P_h=\left(\frac{1}{2}
{\rm I}+K'\right)(V^{-1}-i_h(i_h^* V i_h)^{-1} i_h^*)\left(\frac{1}{2}{\rm I}+K\right). \]
From \cite{MS-1} %(see also \cite{cc}), 
the operator $E_h$ is bounded and there exist a constant $c>0$ such that 
\[
\|E_h(\mathbf{u})\|_{{\mathbf{H}}^{-1/2}(\Gamma)}\leq c \, \displaystyle \inf_{\mathbf{w}\in \mathcal{W}_h}\left \|V^{-1}\left({\frac{1}{2}\rm I}+K\right)\mathbf{u} - w\right\|_{\mathbf{H}^{-1/2}(\Gamma)} \quad \forall \mathbf{u} \in \mathbf{H}^{-1/2}(\Gamma).
\]
The following statements hold:
\begin{lem}\label{lem1}
\begin{description}
\item (i) 
If 
$ \mathbf{u}_h \rightharpoonup \mathbf{u}$ (weak convergence) and $\mathbf{v}_h \to \mathbf{v}$ in $\mathbf{H}^{1/2}(\Gamma)$. Then, $
\displaystyle \lim_{h\to 0} \,\langle P_h \mathbf{u}_h, \mathbf{v}_h \rangle = \langle P \mathbf{u}, \mathbf{v} \rangle.$ %_{\Gamma_0}. $
\item (ii)
If 
$ \mathbf{u}_h \to \mathbf{u}$  and $\mathbf{v}_h \rightharpoonup  \mathbf{v}$ in $\mathbf{H}^{1/2}(\Gamma)$. Then, $
\displaystyle \lim_{h\to 0} \,\langle P_h \mathbf{u}_h, \mathbf{v}_h \rangle = \langle P \mathbf{u}, \mathbf{v} \rangle $. %_{\Gamma_0}. $
\end{description}
\end{lem}
{\bf Proof:} The part $(i)$ follows immediately from the estimate below. Indeed, from \cite[Lemma 9]{cc}, there exists a constant $c_0$ such that
\[
\langle P_h \mathbf{v}_h -i_h^*P \mathbf{v}, \mathbf{w}_h \rangle_{\mathcal{V}_h} \leq c_0 \|\mathbf{w}_h\|_{\mathbf{H}^{1/2}(\Gamma)}\left( e_h(\mathbf{v})+ \|\mathbf{v}_h-\mathbf{v}\|_{\mathbf{H}^{1/2}(\Gamma)}\right) %\quad \forall  \mathbf{v} \in \mathbf{H}^{1/2}(\Gamma), \; \mathbf{v}_h, \mathbf{w}_h \in \mathcal{V}_h,
\]
for any $\mathbf{v} \in \mathbf{H}^{1/2}(\Gamma)$ and for any  $\mathbf{v}_h, \mathbf{w}_h \in \mathcal{V}_h$, where $e_h(\mathbf{v})$ satisfies $e_h(\mathbf{v})\to 0$ as $h\to 0$. 

Hence, using the symmetry  of $P$ and $P_h$, we obtain
\begin{eqnarray*}
\langle P_h \mathbf{u}_h, \mathbf{v}_h \rangle - \langle P \mathbf{u}, \mathbf{v} \rangle & = &\langle P_h \mathbf{v}_h-i_h^*P\mathbf{v}, \mathbf{u}_h \rangle + \langle  P \mathbf{v}, i_h \mathbf{u}_h-\mathbf{u} \rangle \\
& \leq&  c_0 \|\mathbf{u}_h\|_{\mathbf{H}^{1/2}(\Gamma)}\left ( e_h(\mathbf{v}) +\|\mathbf{v}_h-\mathbf{v}\|_{\mathbf{H}^{1/2}(\Gamma)}\right) + \langle  P \mathbf{v},  \mathbf{u}_h-\mathbf{u} \rangle
\end{eqnarray*}
and thus, (i) is satisfies. 
% &= &\langle P \mathbf{u}_h, \mathbf{v}_h \rangle_{\Gamma_0} -  \langle P\mathbf{u}, \mathbf{v} \rangle_{\Gamma_0}  -\langle E_h \mathbf{u}_h, \mathbf{v}_h \rangle_{\Gamma} \\[0.2cm]
%&=& \langle P  \mathbf{u}_h,  \mathbf{v}_h -\mathbf{v} \rangle_{\Gamma_0}+ \langle  P  \mathbf{u}_h -P u,\mathbf{v} \rangle_{\Gamma_0} -\langle E_h \mathbf{u}_h, \mathbf{v}_h \rangle_{\Gamma} \\[0.2cm]
%&\leq & \| P j_h \mathbf{u}_h \|_{\mathbf{H}^{-1/2}(\Gamma)} \| j_h \mathbf{v}_h -\mathbf{v}\|_{\mathbf{H}^{1/2}(\Gamma)} +
%\langle  P j_h \mathbf{u}_h -P u,\mathbf{v} \rangle_{\Gamma_0}\\[0.2cm] 
%&+&\|E_h\|_{\mathcal{L}(\mathbf{H}^{1/2}(\Gamma),\mathbf{H}^{-1/2}(\Gamma))}\|j_h\mathbf{u}_h\|_{\mathbf{H}^{1/2}(\Gamma)}\|j_h\mathbf{v}_h\|_{\mathbf{H}^{1/2}(\Gamma)}. 
%\end{eqnarray*}
%\normalsize
The proof of $(ii)$ follows in the same way. 
\qed

\begin{figure}[t]
\centering
\includegraphics[trim=0cm 12cm 0cm 8cm, clip, width=1.0\textwidth]{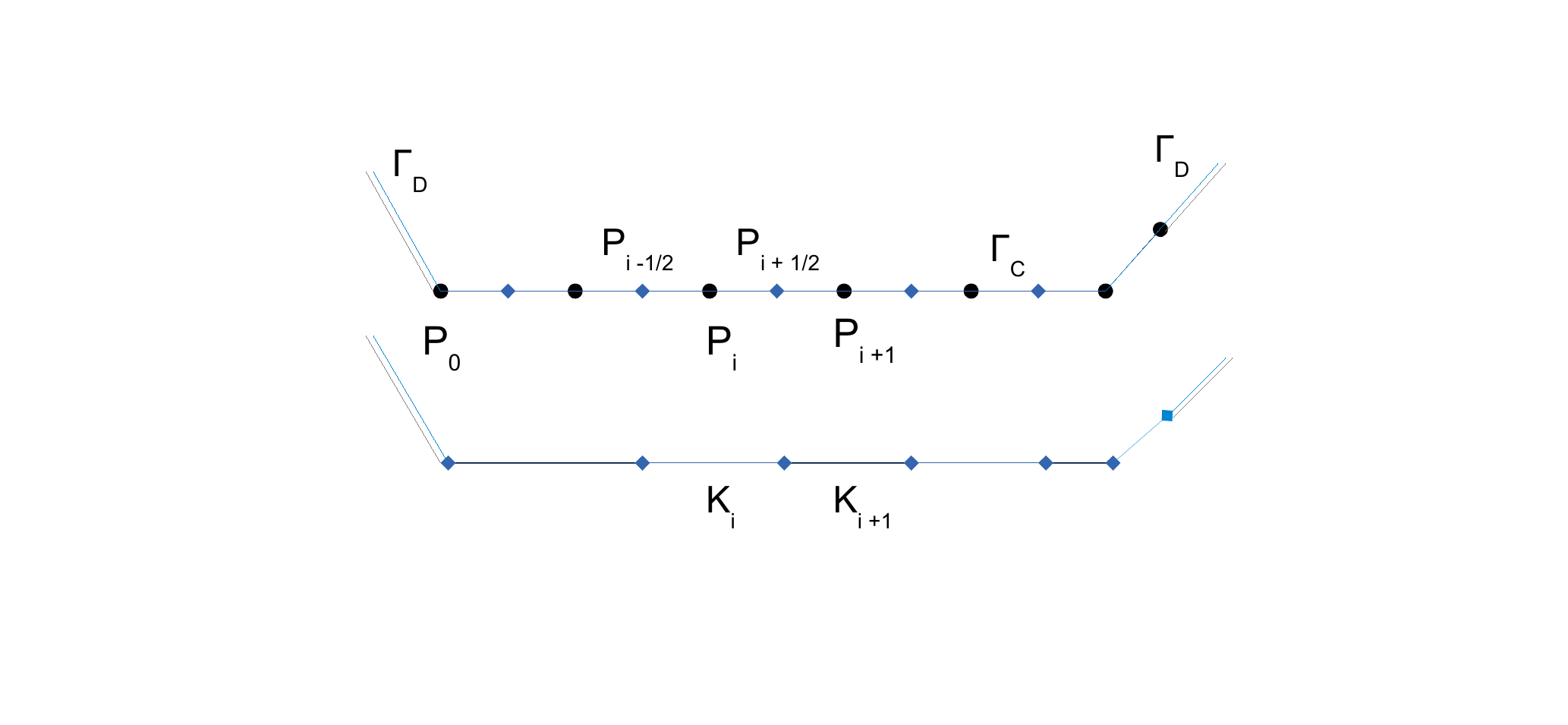}
\vspace{-2cm}
\caption{Discretization on $\Gamma_C$; $\rm{P_0}$ is a boundary point for $\Gamma_D$} 
\label{fig:1}
\end{figure}

From now on, let $\Omega $ be a bounded domain in $\R^2$ with a polygonal boundary $\Gamma$. We define $\Pi : \mathbf{H}^{1/2} (\Gamma) \to L^2(\Gamma_C)$ by 
$$
\Pi \mathbf{u}_h= \mathbf{u}_h \cdot \mathbf{n} \quad \mbox{on} \quad \Gamma_C .
$$
The mapping $\Pi$ is linear continuous from $\mathbf{H}^{1/2} (\Gamma)$ into $L^2(\Gamma_C)$, i.e.
\begin{equation} \label{contPi}
\exists c>0 \, :\, \| \Pi \mathbf{v} \|_{L^2(\Gamma_C)} \leq c \| \mathbf{v}\|_{\mathbf{H}^{1/2} (\Gamma)} \quad \forall v \in \mathbf{H}^{1/2} (\Gamma).
\end{equation}
Further,  let  $\{\Gamma^j\}_{j=1}^J$ be a system of all open straight line segments $\Gamma^j$ of $\Gamma$. We denote by $\tilde{\mathcal{V}}_h$  the image of $\mathcal{V}_h$ with respect to $\Pi$, i.e.
\[
\tilde{\mathcal{V}}_h= \{ w_h \in L^\infty(\overline{\Gamma}_C) \, :\, w_h|_{\Gamma^j}  \in C(\overline{\Gamma^j}) \; \, \forall \Gamma^j, \; \mbox{and} \;  w_h |_{E} \in P_1(E) \; \, \forall E\in \mathcal{T}_h|_{\Gamma_C}, \; w_h=0 \; \mbox{on} \;  \overline{\Gamma}_D\},
\]
where the symbol $\mathcal{T}_h|_{\Gamma_C}$ denotes the partition of $\Gamma_C$ induced by $\mathcal{T}_h$.  Note that since $\Omega$ is polygonal domain, the unit normal vector $\mathbf{n}$ is piecewise constant with a discontinuity at the vertices of $\Omega$. \\[0.1cm]
Let $\{P_i\}_{i=0}^{m}$ by the set of all nodes of $\mathcal{T}_h$ lying on $\overline{\Gamma}_C$. To approximate the  G\^{a}teaux derivative  $\langle D J_\varepsilon (\cdot), \cdot \rangle $ we use a numerical integration due to the Kepler's trapezoidal rule and obtain
\begin{eqnarray*} 
\langle D J_\varepsilon(\mathbf{u}_{h}), \mathbf{v}_{h} \rangle &\approx &
\langle D J_{\varepsilon,h}(\mathbf{u}_{h}), \mathbf{v}_{h} \rangle \nonumber  
 \\ [0.1cm] & = & \frac{1}{2}\sum_{i=0}^{m-1}
|P_iP_{i+1}|\,\big[  S_x(\varepsilon,
   \Pi \mathbf{u}_{h} (P_i)) \Pi
\mathbf{v}_{h}(P_i) % \nonumber \\[0.1cm]
 +   S_x(\varepsilon,
   \Pi \mathbf{u}_{h} (P_{i+1}))  \Pi
\mathbf{v}_{h}(P_{i+1})  \big ] \, .%\label{disc_func}.
\end{eqnarray*}
The discretization of the regularized problem (\ref{reg}) reads now as follows: \\[0.2cm] {\bf Problem} ($\mathcal{P}_{\varepsilon, h}$)~ Find
   $u_{\varepsilon,h} \in K^{\Gamma}_{h}$ such that 
\begin{equation} \label{reg_dis_problem}
\hspace{-1.0cm}
(\mathbf{v}_h-\mathbf{u}_{\varepsilon,h})^T {P}_h \mathbf{u}_{\varepsilon,h} +
\langle D J_{\varepsilon, h}(\mathbf{u}_{\varepsilon,h}), \mathbf{v}_h -\mathbf{u}_{\varepsilon,h} \rangle_{\Gamma_C} 
 \geq \int_{\Gamma_0} N\mathbf{f} \cdot (\mathbf{v}_h-\mathbf{u}_{\varepsilon,h}) ds+\int_{\Gamma_N} \mathbf{t} \cdot (\mathbf{v}_h-\mathbf{u}_{\varepsilon, h}) \, ds, \quad \forall \mathbf{v}_h \in \mathcal{K}_h^{\Gamma}.  
\end{equation}
\normalsize
Let $\mathcal{D}_h$ be another partition of $\Gamma_C$  consisting of elements $K_i$ joining the midpoints $P_{i-1/2}$, $P_{i+1/2}$ of the edges $E \in \mathcal{T}_h$ lying on $\Gamma_C$ sharing $P_{i}$ as a common point. If $P_i$ is a vertex of $\partial \Omega$ then $K_i$ is the half of the edge. Moreover, if the segment $K_i$ is adjacent to the boundary node $P_i$ of $\overline{\Gamma}_D$, it will be appended to its neighbour $K_{i+1}$, see Figure  \ref{fig:1}. 
Further, on $\mathcal{D}_h$ we introduce the space $\mathcal{Y}_h$ of all piecewise constant functions by 
\[
\mathcal{Y}_h=\{\mu_h \in L^{\infty}(\Gamma_C)\, :\, \mu_h|_K \in P_0(K) \quad \forall K \in \mathcal{D}_h\}
\]
and define the piecewise constant Lagrange interpolation operator $
L_h \,:\, \tilde{\mathcal{V}}_h \to \mathcal{Y}_h$ by
\[
L_h(w_h)(x)=\sum _{i} %{P_i \in \Gamma_C \cap \Sigma_h}
 w_h(P_i) \,
\chi_{\mbox{\small int} \;_{\Gamma_C} K_i}(x),
\]
where $\chi_{\small \mbox{int} \;_{\Gamma_C} K_i}$ is the characteristic function of the interior of $K_i$ in $\Gamma_C$.
 
It holds that
\begin{equation} \label{exint}
\langle D J_{\varepsilon,h}(\mathbf{u}_{h}), \mathbf{v}_{h} \rangle_{\Gamma_C} = \displaystyle  {\int_{\Gamma_C}}  
S_x(\varepsilon, L_h( \Pi \mathbf{u}_{h})) L_h(\Pi \mathbf{v}_{h})  \,ds. 
\end{equation}
%Moreover, the operator $D J_{\varepsilon,h} \,:\,  \mathcal{V}_h \to \mathcal{V}^*_h$ is strongly continuous, and 
By (\ref{ass_4}), there exists a constant $C>0$ independent of $\varepsilon$ and $h$ such that
\begin{equation}\label{est01}
\exists C >0 \, :\, \langle D J_{\varepsilon,h}(\mathbf{u}_{h}),\mathbf{u}_{h}\rangle_{\Gamma_C} \geq -C \|\mathbf{u}_{h}\|_{\mathbf{{H}}^{1/2}(\Gamma)}.
\end{equation}
%for the proofs see \cite{Ovcharova2012}.  

From Glowinski et al. \cite{Gl_1981}, we know that
\begin{equation}%(lcl)
 \|L_h ( \mathbf{v}_h \cdot \mathbf{n})\|_{\mathbf{L}^2(\Gamma)}  \leq   2\,\|  \mathbf{v}_{h} \cdot \mathbf{n}\|_{\mathbf{L}^2(\Gamma)}  \label{Gl1_1} %\\[0.3cm]
\end{equation}
and therefore, 
\begin{equation}%(lcl)
 \|\mathbf{v}_h \cdot \mathbf{n} - L_h ( \mathbf{v}_h \cdot \mathbf{n})\|_{\mathbf{L}^2(\Gamma)}  \leq   3\,\|  \mathbf{v}_{h} \cdot \mathbf{n}\|_{\mathbf{L}^2(\Gamma)}.   \label{Gl1} %\\[0.3cm]
\end{equation}
%{\color{blue}
Let now $\mathbf{H}^s(\Gamma^j)$, $s\geq 0$, be the standard Sobolev space from \cite{Grisvard,Grisvard2011,Lions} defined on the open straight pieces $\Gamma^j$  by 
\[
\mathbf{H}^s(\Gamma^j)=\{u|_{\Gamma^j}\; :\; u \in \mathbf{H}^s(\Gamma)\}.
\]
According to Grisvard \cite{Grisvard,Grisvard2011}, $\mathbf{H}^s(\Gamma)\subset \displaystyle \prod _{j=1}^J \mathbf{H}^s(\Gamma^j)$ for $s\in [1/2, 3/2)$ and 
\begin{equation} \label{Grisvard1}
\displaystyle \sum _{j=1}^J \|\mathbf{u}\|^2_{\mathbf{H}^s(\Gamma^j)} \leq C \| \mathbf{u} \|^2_{\mathbf{H}^s(\Gamma)}.
\end{equation}
Again from  Glowinski et al. \cite{Gl_1981},
\begin{equation}\label{Gl3}
\|\mathbf{v}_h \cdot \mathbf{n} - L_h(\mathbf{v}_h \cdot \mathbf{n})\|^2_{L^2(\Gamma^j)}\leq C h^2 \|\mathbf{v}_h\cdot \mathbf{n}\|^2_{H^1(\Gamma^j)}\leq C h^2 \|\mathbf{v}_h\|^2_{H^1(\Gamma^j)}.
\end{equation}
Summing over all $j$ and using thereafter %, and combining (\ref{Grisvard1}) and 
(\ref{Grisvard1}), it follows that
\begin{equation}\label{Gln4}
\|\mathbf{v}_h \cdot \mathbf{n} - L_h(\mathbf{v}_h \cdot \mathbf{n})\|_{L^2(\Gamma)}\leq C h \|\mathbf{v}_h\|_{H^1(\Gamma)}.
\end{equation}
By interpolation between $L^2(\Gamma)$ and $H^1(\Gamma)$ we deduce from (\ref{Gl1}) and (\ref{Gln4}) that
\begin{equation}\label{Gl4}
\|\mathbf{v}_h \cdot \mathbf{n} - L_h(\mathbf{v}_h \cdot \mathbf{n})\|_{L^2(\Gamma)}\leq C h^{1/2} \|\mathbf{v}_h\|_{H^{1/2}(\Gamma)}.
\end{equation}
By the compactness of $H^{s_1}(\Gamma)\subset H^{s_2}(\Gamma)$ for $0\leq s_2\leq s_1$ ($\Omega \subset \R^2$), this gives
\begin{equation}
\mathbf{v}_h \rightharpoonup \mathbf{v} \; \mbox{in} \; \mathbf{H}^{1/2}(\Gamma) 
\Rightarrow \|L_h( \mathbf{v}_h \cdot \mathbf{n}) - \mathbf{v} \cdot \mathbf{n}\|_{\mathbf{L}^2(\Gamma)} \to  0. \label{Gl3}
\end{equation}
For the proof of the compact embedding $H^{s_1}(\Gamma)\subset H^{s_2}(\Gamma)$ by Fourier expansion, see Kress \cite{Kress}.

%Further, we introduce the functional $\varphi_{\varepsilon,h} \, :\, \mathcal{V}_h \times \mathcal{V}_h \to \R$ by 
%\begin{equation} \label{pfi_descret}
%\varphi_{\varepsilon,h}(\mathbf{u}_{h},\mathbf{v}_{h}) = \langle D J_{\varepsilon,h}(\mathbf{u}_{h}),\mathbf{v}_{h}-\mathbf{u}_{h}\rangle_{\Gamma_C}.
%\end{equation}
%According to \cite{Ovcharova2012}, this functional is pseudomonotone. % and upper semicontinuous with respect to the first argument.
The solvability of $(\mathcal{P}_{\varepsilon,h})$ and the convergence of its solutions to a solution of the boundary hemivariational inequality (\ref{BdHemVar}) relies on the following {\it general approximation result} from  \cite{Gw_Ov}.

Let $\mathcal{K}$ be a closed convex nonvoid subset of a reflexive Banach space $X$. We consider the variational inequality $VI(\psi,f,\mathcal{K})$:\; Find $u \in \mathcal{K}$ such that
\[
\psi(u,v) \geq \langle f, v-u\rangle \quad \forall v \in \mathcal{K}. 
\]
Let $T$ be a directed set.  We introduce the family $\{\mathcal{K}_t\}_{t\in T}$ of nonempty, closed and convex sets $\mathcal{K}_t$ (not necessary contained in $\mathcal{K}$) and assume the following hypotheses:
\begin{description}
\item (H1) If  $\{v_{t'}\}_{t'\in T'}$ weakly converges to $v$ in $X$, $v_{t'} \in
  \mathcal{K}_{t'} ~(t' \in T') $ for a subnet $\{\mathcal{K}_{t'}\}_{t'\in T'}$ of the  net  $\{\mathcal{K}_t\}_{t\in T}$,  then   $v\in \mathcal{K}$.
\item (H2) For any $v\in \mathcal{K}$ and any $t \in T$ there exists $v_{t} \in \mathcal{K}_{t}$ 
such that  $v_{t} \to v$  in $X$.
\item (H3) $\psi_t$ is pseudomonotone for any $t \in T$.
\item (H4) $f_t \to f$ in $V^\ast$.
\item (H5) For any nets $\{u_t\}$ and $\{v_t\}$ 
 such that  $u_t \in \mathcal{K}_t$,   $v_t\in \mathcal{K}_t$, 
  $u_t \rightharpoonup u$, and $v_t \to v$ in $X$ it follows that 
$$
\displaystyle {\limsup_{t \in T}}\,~ \psi_t(u_t,v_t) \leq \psi(u,v) \,.
$$
\item  (H6) 
  There exist constants $c>0$,   $d$, $d_0 \in \R$ and $\alpha >1$  (independent of   $t \in T$) such that for some $w_t \in \mathcal{K}_t$
with $w_t \to w$ there holds  
$$
-\psi_t(u_t,w_t)\geq c \|u_t\|^\alpha_V +d \|u_t\|_V +d_0, \quad \forall u_t \in \mathcal{K}_t, \; \forall t \in T \,.
$$
\end{description}
\begin{theo} \cite{Ov_Gw} \label{theo1}({\it General approximation result}) Under the hypotheses $(H1)$-$(H6)$, there exist a solution $u_t$ to the approximate problem $VI(\psi_t,f_t,\mathcal{K}_t)$ and  the family $\{u_t\}$ is   bounded in $X$. Moreover, there exists a subnet of $\{u_t\}$ that converges weakly  in $X$ to  a solution of  the problem $VI(\psi,f,\mathcal{K})$. Furthermore, any weak accumulation point of  $\{u_t\}$ is a solution to the problem $VI(\psi,f,\mathcal{K})$.
\end{theo}
%{\color{blue}
\begin{rem} 
The hypotheses (H1) and (H2) are due to Glowinski \cite{Gl} and describe the Mosco convergence \cite{At} of the family $\mathcal{K}_t$ to $\mathcal{K}$. If we replace the weak convergence in (H1) by the strong convergence, we obtain
the classical Painlev\'{e}-Kuratowski convergence \cite{Aubin} based on the notions of upper and lower limits of a net of sets. As already mentioned, the pseudomonotonicity of $\psi_t$ in (H3) guarantees a solution to the discrete problem, whereas $f_t$ in (H4) is a standard approximation of the linear functional $f$, for example, by  numerical integration. The verification of (H5) and (H6) %for the functionals $\varphi_{\varepsilon,h}$ defined by (\ref{pfi_descret}) 
is based on the arguments used in \cite{Ovcharova2012}, where the  discretization of the domain hemivariational problem has been investigated. For self-consistency of the paper we include them below. 
\end{rem}

{\bf Verification of $(H5)$} 
Let 
\[\mathbf{u}_{\varepsilon_n,h_n}\rightharpoonup \mathbf{u}, \quad  \mathbf{v}_{\varepsilon_n,h_n}\to \mathbf{v} \quad \mbox{in} \; \mathbf{H}^{1/2}(\Gamma)
\]
 with $h_n\to 0$ and $\varepsilon_n\to 0^+$. 
We define
\[
w_{\varepsilon_n,h_n}=L_{h_n}(\Pi (\mathbf{v}_{\varepsilon_n,h_n}-\mathbf{u}_{\varepsilon_n,h_n})). %\quad \mbox{and} \quad w=\gamma (v-u)= \gamma v- \gamma u.
\]
Since $\mathbf{v}_{\varepsilon_n,h_n}-\mathbf{u}_{\varepsilon_n,h_n}\rightharpoonup \mathbf{v}-\mathbf{u}=:\mathbf{w}$ in $\mathbf{H}^{1/2}(\Gamma)$, it follows by (\ref{Gl3}) that  
\begin{equation} \label{H5_1}
w_{\varepsilon_n,h_n} \to \mathbf{w} \cdot \mathbf{n} = \Pi \mathbf{w}\quad \mbox{in}
\;L^2(\Gamma) \quad \mbox{as} \; n \to \infty.
\end{equation}
By applying the Cauchy-Schwarz inequality we obtain
\begin{align}
   \int_{\Gamma_C} S_x(L_{h_n}(\varepsilon_n, \Pi \mathbf{u}_{\varepsilon_n,h_n}))
  w_{\varepsilon_n,h_n}    \ ds 
 & =    \int_{\Gamma_C} S_x(\varepsilon_n, L_{h_n}(\Pi \mathbf{u}_{\varepsilon_n,h_n}))   (w_{\varepsilon_n,h_n}-\Pi \mathbf{w})  \, ds   %\nonumber\\
+ \int_{\Gamma_C} S_x(\varepsilon_n, L_{h_n}(\Pi
    \mathbf{u}_{\varepsilon_n,h_n})) \, \Pi \mathbf{w} \, ds  \nonumber\\
&\leq  \| S_x(\varepsilon_n, L_{h_n}(\Pi \mathbf{u}_{\varepsilon_n,h_n}))\|_{L^2(\Gamma_C)}\|w_{\varepsilon_n,h_n}-\Pi \mathbf{w}\|_{L^2(\Gamma_C)}
  %\nonumber\\ 
+    \int_{\Gamma_C}  S_x(\varepsilon_n, L_{h_n}(\Pi \mathbf{u}_{\varepsilon_n,h_n})) \, \Pi \mathbf{w} \,
   ds. \label{H5_2}
\end{align}
By (\ref{Gl3}), $L_{h_n}(\Pi \mathbf{u}_{\varepsilon_n,h_n}) \to \Pi \mathbf{u}$ in $L^2(\Gamma_C)$ and thus, for a subsequence, 
$
L_{h_n}(\Pi \mathbf{u}_{\varepsilon_n,h_n})(s) \to \Pi \mathbf{u} (s) $  for a.e.  $s\in \Gamma_C $. Moreover, there exists a function $h\in L^2(\Gamma_C;\R_+)$ such that $|L_{h_n}(\Pi \mathbf{u}_{\varepsilon_n,h_n})(s)|\leq h(s)$ for a.e. $s\in \Gamma_C$. % with $h\in L^2(\Gamma_C;\R_+)$. 

The last integrand in (\ref{H5_2}) is uniformly bounded from above, since by (\ref{ass_3}),
\[
 | S_x(L_{h_n}(\varepsilon_n,\Pi \mathbf{u}_{\varepsilon_n,h_n})(s)) \, \Pi \mathbf{w}(s)  | \leq c(1+|L_{h_n}(\Pi \mathbf{u}_{\varepsilon_n,h_n})(s)|) |\Pi \mathbf{w}(s)|\leq c(1+h(s))|\Pi \mathbf{w}(s)|.
\]
Moreover, applying $(a+b)^2\leq 2(a^2+b^2)$ to the right-hand side of (\ref{ass_3}) and integrating over $\Gamma_C$ implies
\begin{equation} \label{app1}
\int_ {\Gamma_C} | S_x( \varepsilon_n,  L_{h_n}(\varepsilon_n,\Pi \mathbf{u}_{\varepsilon_n,h_n}))|^2 \, ds \leq 2c^2
\,\mbox{meas}\, (\Gamma_C) + 2c^2\,\| L_{h_n}(\varepsilon_n,\Pi \mathbf{u}_{\varepsilon_n,h_n}) \|^2_{L^2(\Gamma_C)} \leq \tilde{c} (1+\|h\|^2_{L^2(\Gamma_C)}).
\end{equation}
We are now in position to apply  Fatou's lemma and in view of (\ref{ass_3b}) we get
\begin{equation}\label{H5_5}
\limsup_{n\to \infty}\int_{\Gamma_C}  S_x(\varepsilon_n, L_{h_n}(\Pi \mathbf{u}_{\varepsilon_n,h_n})) \, \Pi \mathbf{w} \,
   ds \leq \int _{\Gamma_C} \limsup_{n\to \infty} S_x(\varepsilon_n, L_{h_n}(\Pi \mathbf{u}_{\varepsilon_n,h_n})) \, \Pi \mathbf{w} \,
   ds\leq \int_{\Gamma_C} f^0(\Pi \mathbf{u} (s);\Pi \mathbf{w}(s))\, ds. %\quad \forall \tilde{v}\in L^2(\Gamma_C).
\end{equation}
Passing to limsup in (\ref{H5_2}) we finally conclude that   
 \[
  \displaystyle {\limsup_{n\to \infty}} \,\int_{\Gamma_C} S_x (\varepsilon_n, L_{h_n}(\Pi  \mathbf{u}_{\varepsilon_n,h_n})) \, L_{h_n} (\Pi \mathbf{v}_{\varepsilon_n,h_n}- \Pi \mathbf{u}_{\varepsilon_n,h_n})\, ds \\
  \leq    \int_{\Gamma_C} f^0(\Pi \mathbf{u} ;\Pi \mathbf{v}- \Pi \mathbf{u}) \, ds,
\]
where we have used (\ref{H5_5}),  the strong convergence (\ref{H5_1}) as well as the boundedness of $\{S_x (L_{h_n}(\Pi \mathbf{u}_{\varepsilon_n,h_n}),\varepsilon_n)\}$ in $L^2(\Gamma_C)$ (see (\ref{app1})). Thus, the hypothesis (H5) is verified for $t=(\varepsilon_n, h_n)$. \qed %and $\psi_t:=\varphi_{\varepsilon_n, h_n}$.\qed
\\
{\bf Verification of $(H6)$} The hypothesis (H6) is obviously 
satisfied, since  by (\ref{est01}) there exists a  constant $C>0$, which does not depend on  $\varepsilon$ and $h$, such that
\[
%\varphi_{\varepsilon_n,h_n}(\mathbf{u}_{\varepsilon_n,h_n},0)= 
\langle D J_{\varepsilon_n,h_n} (\mathbf{u}_{\varepsilon_n,h_n}), -\mathbf{u}_{\varepsilon_n,h_n}
\rangle_{\Gamma_C} \leq C\|\mathbf{u}_{\varepsilon_n,h_n}\|_{\mathbf{H}^{1/2}(\Gamma)}.
\]
For the  convenience of the reader, we next show the uniform boundedness of $\{\mathbf{u}_{\varepsilon,h}\}$ in $\mathbf{\tilde{H}}^{1/2}(\Gamma_0)$. 
\begin{lem} \label{uniform}
The family $\{\mathbf{u}_{\varepsilon,h}\}$ of solutions of the problem $(P_{\varepsilon, h})$ is uniformly bounded in $\mathbf{\tilde{H}}^{1/2}(\Gamma_0)$. 
\end{lem}
{\bf Proof}\, 
Putting $\mathbf{v}_{h}=0$ in (\ref{reg_dis_problem}), using  (\ref{coerc})  
and  the estimate (\ref{est01}), we get
%\begin{eqnarray*}
\[
c\|\mathbf{u}_{\varepsilon,h}\|^2_{\mathcal{V}}   \leq \langle P_h \mathbf{u}_{\varepsilon,h}, \mathbf{u}_{\varepsilon,h} \rangle_{\Gamma_0} \leq %\langle \mathbf{g},\mathbf{u}_{\varepsilon,h}\rangle + \varphi_{\varepsilon,h}(\mathbf{u}_{\varepsilon,h},0) \\[0.2cm]
\langle \mathbf{g},\mathbf{u}_{\varepsilon,h}\rangle + \langle D J_{\varepsilon,h}(\mathbf{u}_{\varepsilon,h}), -\mathbf{u}_{\varepsilon,h} \rangle_{\Gamma_C} 
\leq  \|\mathbf{g}\|_{\mathcal{V}^*} \,\|\mathbf{u}_{\varepsilon,h}\|_{\mathcal V} + c \|\mathbf{u}_{\varepsilon,h}\|_{\mathcal V}, 
%\end{eqnarray*}
\]
which  implies the uniform boundedness of $\{\mathbf{u}_{\varepsilon,h}\} $ in $\mathcal{V}$ with respect to the both parameters $\varepsilon$ and $h$.
\qed
\\[0.2cm]
Further, in case of uniqueness we improve the convergence result of Theorem \ref{theo1} and show that the weak convergence can be replaced by the strong one. 
\begin{theo} \label{lem2}
Let the solutions $\mathbf{u}$ to $(\mathcal{P})$ and $\mathbf{u}_{\varepsilon, h}$ to $(\mathcal{P}_{\varepsilon,h})$ exist uniquely. Then 
\[
\displaystyle \lim _{\varepsilon \to 0, h \to 0} \| \mathbf{u}_{\varepsilon, h} - \mathbf{u} \|_{\tilde{\mathbf{H}}^{1/2}(\Gamma_0)} =0.
\]
\end{theo}
{\bf Proof of Theorem \ref{lem2}}
Let $\{h_n\}$  and $\{\varepsilon_n\}$ be  arbitrary sequences such that $h_n\to 0^+$ and $\varepsilon_n \to 0 ^+$ as $n\to \infty$. In view of $\mathbf{(ii)}$,  there exists a sequence $\{\bar{\mathbf{u}}_{\varepsilon_n,h_n}\}$ such that 
$\bar{\mathbf{u}}_{\varepsilon_n,h_n} \in \mathcal{K}^{\Gamma}_{h_n}$ and $\bar{\mathbf{u}}_{\varepsilon_n,h_n} \to \mathbf{u}$ in $\mathcal{V}$. 

Using (\ref{coerc}), %(here we have dropped the symbol $j_h$ on the left-hand side), 
we obtain 
\begin{eqnarray} \hspace{-1cm}
c\|\bar{\mathbf{u}}_{\varepsilon_n,h_n}-\mathbf{u}_{\varepsilon_n,h_n}\|^2_{\mathcal{V}} & \leq & \langle P_h (\bar{\mathbf{u}}_{\varepsilon_n,h_n}-\mathbf{u}_{\varepsilon_n,h_n}),
\bar{\mathbf{u}}_{\varepsilon_n,h_n}-\mathbf{u}_{\varepsilon_n,h_n} \rangle_{\Gamma_0} \nonumber \\[0.2cm]
& = &\langle  P_h \bar{\mathbf{u}}_{\varepsilon_n,h_n}, \bar{\mathbf{u}}_{\varepsilon_n,h_n}-\mathbf{u}_{\varepsilon_n,h_n} \rangle_{\Gamma_0} - \langle  P_h \mathbf{u}_{\varepsilon_n,h_n},
\bar{\mathbf{u}}_{\varepsilon_n,h_n}-\mathbf{u}_{\varepsilon_n,h_n} \rangle_{\Gamma_0}. \label{strongest}
\end{eqnarray}
Since $\bar{\mathbf{u}}_{\varepsilon_n,h_n} \to \mathbf{u}$ in $\mathcal{V}$ and ${\mathbf{u}}_{\varepsilon_n,h_n} \rightharpoonup \mathbf{u}$ in $\mathcal{V}$, it follows from Lemma \ref{lem1} $(ii)$ that the first term on the right-hand side of (\ref{strongest}) tends to zero.

Using the definition of $(P_{\varepsilon_n,h_n})$, inequality (\ref{reg_dis_problem}), the second term 
can be estimated as follows:
\begin{equation} \label{est002}
| \langle  P_h \mathbf{u}_{\varepsilon_n,h_n}, \mathbf{u}_{\varepsilon_n,h_n} - \bar{\mathbf{u}}_{\varepsilon_n,h_n} \rangle_{\Gamma_0} | 
\leq  
 |\langle \mathbf{g}, \mathbf{u}_{\varepsilon_n,h_n} - \bar{\mathbf{u}}_{\varepsilon_n,h_n} \rangle | + |\langle D J_{\varepsilon_n,h_n}(\mathbf{u}_{\varepsilon_n,h_n}), \bar{\mathbf{u}}_{\varepsilon_n,h_n}-\mathbf{u}_{\varepsilon_n,h_n} \rangle_{\Gamma_C} |.
 \end{equation}
In addition, we have
\begin{eqnarray*}
 \hspace{-1.5cm} |\langle D J_{\varepsilon_n,h_n}(\mathbf{u}_{\varepsilon_n,h_n}),
\bar{\mathbf{u}}_{n,h_n}-\mathbf{u}_{n,h_n} \rangle_{\Gamma_C} | 
&= &
\left | {\int_{\Gamma_C}}  
S_x( L_{h_n}(\varepsilon_n, \Pi \mathbf{u}_{\varepsilon_n,h_n})) \,  L_{h_n}(\Pi (\bar{\mathbf{u}}_{\varepsilon_n,h_n}- \mathbf{u}_{\varepsilon_n,h_n}))  \,ds \right | \\ [0.2cm] & 				\leq &
\| 
S_x(\varepsilon_n, L_{h_n}(\Pi \mathbf{u}_{\varepsilon_n,h_n}))\| _{L^2(\Gamma_C)} \| L_{h_n}(\Pi (
\bar{\mathbf{u}}_{\varepsilon_n,h_n}- \mathbf{u}_{\varepsilon_n,h_n}))\| _{L^2(\Gamma_C)} \to 0,
\end{eqnarray*}
as follows from (\ref{app1}), the boundedness of $\{\mathbf{u}_{\varepsilon_n,h_n}\}$ in $\tilde{\mathbf{H}}^{1/2}(\Gamma_0)$ and (\ref{Gl3}). 

Passing now to the limit superior in (\ref{est002}), we get
$$
\limsup_{n\to \infty} \, \langle  P_h  \mathbf{u}_{\varepsilon_n,h_n},\mathbf{u}_{\varepsilon_n,h_n}- \bar{\mathbf{u}}_{\varepsilon_n,h_n} \rangle _{\Gamma_0} \leq 0.
$$
Hence, (\ref{strongest}) entails in the limit, that
$$
\limsup_{n\to \infty}  c \|\bar{\mathbf{u}}_{\varepsilon_n,h_n}- \mathbf{u}_{\varepsilon_n,h_n }\|^2_{\mathcal{V}}  \leq 0
$$
and therefore, 
$$ \|\bar{\mathbf{u}}_{\varepsilon_n,h_n} -\mathbf{u}_{\varepsilon_n,h_n}\|_\mathcal{V}\to 0.$$
Finally, from  the triangle inequality 
$$
\|\mathbf{u}_{\varepsilon_n,h_n}-\mathbf{u} \|_{\mathcal{V}}\leq \|\mathbf{u}_{\varepsilon_n,h_n}- \bar{\mathbf{u}}_{\varepsilon_n,h_n}\|_\mathcal{V} + \|  \bar{\mathbf{u}}_{\varepsilon_n,h_n}-u\| _\mathcal{V},
$$
we get the strong convergence of  an appropriate subsequence of $\{\mathbf{u}_{\varepsilon_n,h_n}\}$ to $\mathbf{u}$ in $\mathcal{V}$. \qed

\section{A-priori error estimate}
%{\color{blue}
In this section we  present an abstract C\'{e}a-Falk approximation lemma for the regularized problem. %, which splits the error under study into four different error terms.  
For its proof we slightly extend the arguments of Maischak and Stephan in \cite {MS-1} for Signorini contact to include the approximation of $D J_\varepsilon$ by $D J_{\varepsilon,h}$. We apply this lemma to obtain an a-priori error estimate for the $h$ - approximate solution of the regularized problem assuming $\mathbf{H}^{3/2}(\Gamma)$ regularity of the solution $\mathbf{u}_\varepsilon$. For our more general problem we arrive at the same convergence rate of $\mathcal{O}(h^{1/4})$ as in \cite {MS-1}.

In addition we refer to Eck et al. \cite{Eck}. They obtain a sharper error estimate for the approximation of the regularized solution of the Coulomb friction problem. However, one should note that the treatment of the  Coulomb friction involves the regularization of the absolute value function only, whereas in the delamination problem we have to cope with multivalued laws and several jumps. 

\begin{lem} \label{Cea} Let $\mathbf{u}_\varepsilon \in \mathcal{K}^\Gamma$,  $\mathbf{u}_{\varepsilon,h} \in \mathcal{K}_{h}^\Gamma$ be  the solutions of the problems $(\mathcal{P}_\varepsilon)$ and  $(\mathcal{P}_{\varepsilon, h})$, respectively.
%Let $\mathbf{u}_\varepsilon \in \mathcal{K}^\Gamma$ be the  solution of the problem $(\mathcal{P}_\varepsilon)$ and let  $\mathbf{u}_{\varepsilon,h} \in \mathcal{K}_{h}^\Gamma$ be the solution of the problem $(\mathcal{P}_{\varepsilon, h})$. 
Assume that %$\mathbf{u}_\varepsilon$ $\in \mathbf{H}^{3/2}(\Gamma)$,  and 
$P \mathbf{u}_\varepsilon -\mathbf{g} \in L^2(\Gamma)$.  Then,  there exists a positive constant $C$ independent of $\varepsilon$ and $h$ such that
\begin{align}
c_P\|\mathbf{u}_\varepsilon-\mathbf{u}_{\varepsilon, h}\|^2_{\mathbf{H}^{1/2}(\Gamma)}  & \leq C \big\{ \|E_{h}(\mathbf{u}_\varepsilon)\|^2_{\mathbf{H}^{-1/2}(\Gamma)}
 + %\inf _{\mathbf{v}\in \mathcal{K}^\Gamma} \{
\|P \mathbf{u}_\varepsilon-\mathbf{g}\|_{\mathbf{L}^2(\Gamma)} \|\mathbf{u}_{\varepsilon,h}-\mathbf{v}\|_{\mathbf{L}^2(\Gamma)} 
 \nonumber \\
& +   %\inf _{\mathbf{v}_h\in \mathcal{K}_h^\Gamma} \{
\|\mathbf{u}_\varepsilon- \mathbf{v}_{ h}\|_{\mathbf{H}^{1/2}(\Gamma)}^2 + \|P \mathbf{u}_\varepsilon-\mathbf{g}\|_{\mathbf{L}^2(\Gamma)} \|\mathbf{u}_\varepsilon-\mathbf{v}_{ h}\|_{\mathbf{L}^2(\Gamma)} \nonumber\\
&+ \langle D J_\varepsilon(\mathbf{u}_\varepsilon),\mathbf{v}-\mathbf{u}_{\varepsilon} \rangle_{\Gamma_C} + \langle D J_{\varepsilon,h}(\mathbf{u}_{\varepsilon, h}),\mathbf{v}_{h}-\mathbf{u}_{\varepsilon,h} \rangle_{\Gamma_C}  \big\}
\label{est1}
\end{align}
for all $\mathbf{v}\in \mathcal{K}^\Gamma$ and for all $\mathbf{v}_h\in \mathcal{K}_h^\Gamma$.
\end{lem}
{\bf Proof} The proof follows by the definitions of the problems $(\mathcal{P}_\varepsilon)$ and $(\mathcal{P}_{\varepsilon,h})$, and by using  estimates similar to (25)-(28) in \cite[Theorem 3]{MS-1}. \qed
\begin{theo} %Suppose (\ref{uniq_reg}).
 Let $\mathbf{u}_\varepsilon \in \mathcal{K}^\Gamma$,  $\mathbf{u}_{\varepsilon,h} \in \mathcal{K}_{h}^\Gamma$ be  the solutions of the problems $(\mathcal{P}_\varepsilon)$ and  $(\mathcal{P}_{\varepsilon, h})$, respectively.  Assume that %$\alpha<c_p/2$, 
$\mathbf{u}_\varepsilon$ $\in \mathbf{H}^{3/2}(\Gamma)$ and $P \mathbf{u}_\varepsilon -\mathbf{g} \in \mathbf{L}^2(\Gamma)$.  Then, under the assumption (\ref{uniq_reg}), there exists a constant  $c=c(\mathbf{u}_\varepsilon, \mathbf{g})$  %\mathbf{t})>0$ 
independent of $h$ such that
%\begin{equation} \label{a-priori_est}
\[
\|\mathbf{u}_\varepsilon-\mathbf{u}_{\varepsilon,h}\|_{\mathbf{H}^{1/2}(\Gamma)} \leq c h^{1/4}.
\]
%\end{equation}
\end{theo} %\cite{Gwinner_2013, MS-1}
{\bf Proof} %To estimate the first term in (\ref{est1}),
We apply Lemma \ref{Cea} with $\mathbf{v}=\mathbf{u}_{\varepsilon,h} \in \mathcal{K}_h^\Gamma \subset \mathcal{K}^\Gamma$ and $\mathbf{v}_h=i_{h} \mathbf{u}_\varepsilon \in \mathcal{K}_h^\Gamma$, the piecewise linear interpolate of $\mathbf{u}_\varepsilon\in \mathbf{H}^{3/2}(\Gamma)\subset \mathbf{C}^0(\Gamma)$. 
According to \cite[Lemma 5]{MS-1}, the first term in (\ref{est1}) can be estimated by %there holds
\begin{equation} \label{Eh}
\|E_h \mathbf{u}_\varepsilon\|_{\mathbf{H}^{-1/2}(\Gamma)} \leq c h \|\mathbf{u}_\varepsilon\|_{\mathbf{H}^{3/2}(\Gamma)}.
\end{equation}
Further, by \cite{MS-1}, proof of Theorem 3, there exists a constant $c>0$ independent of $h$ such that
\begin{equation} \label{est_apriori1}
 \|\mathbf{u}_\varepsilon-i_h \mathbf{u}_\varepsilon\| _{H^{1/2}(\Gamma)} \leq c h \| \mathbf{u}_\varepsilon \|_{\mathbf{H}^{3/2}(\Gamma)}.
 \end{equation}
Since the consistency error $\inf \{ \ldots  \,| \, \mathbf{v} \in \mathcal{K}^\Gamma\}$ disappears, to complete the proof it remains to estimate the last error term in (\ref{est1}). For this purpose, using the Cauchy-Schwarz inequality and assumption (\ref{uniq_reg}) guaranteeing the uniqueness of the solution, we proceed in the following way:
\begin{align*} 
D:=\langle &  D J_{\varepsilon}(\mathbf{u}_{\varepsilon}), \mathbf{u}_{\varepsilon, h}-\mathbf{u}_{\varepsilon} \rangle_{\Gamma_C} + \langle D J_{\varepsilon, h}(\mathbf{u}_{\varepsilon, h}),i_{h} \mathbf{u}_\varepsilon-\mathbf{u}_{\varepsilon} \rangle_{\Gamma_C}  \nonumber \\ & = \int_{\Gamma_C} S_x(\varepsilon,\Pi \mathbf{u}_{\varepsilon})(\Pi \mathbf{u}_{\varepsilon, h}- \Pi \mathbf{u}_{\varepsilon}) \, ds +  \int_{\Gamma_C} S_x(\varepsilon, L_h(\Pi \mathbf{u}_{\varepsilon, h}))L_h(\Pi i_h \mathbf{u}_{\varepsilon}- \Pi \mathbf{u}_{\varepsilon}) \, ds \\
& = \int_{\Gamma_C} S_x(\varepsilon, \Pi \mathbf{u}_{\varepsilon})(\Pi \mathbf{u}_{\varepsilon, h}- L_h(\Pi \mathbf{u}_{\varepsilon, h})) \, ds  +  \int_{\Gamma_C} S_x (\varepsilon, L_h(\Pi \mathbf{u}_{\varepsilon, h})) (L_h (\Pi i_h \mathbf{u}_{\varepsilon})- \Pi \mathbf{u}_{\varepsilon}) \, ds \\
& + \int_{\Gamma_C} (S_x( \varepsilon, \Pi \mathbf{u}_{\varepsilon})-S_x (\varepsilon, L_h(\Pi \mathbf{u}_{\varepsilon, h})) )\, (L_h(\Pi \mathbf{u}_{\varepsilon, h})-\Pi \mathbf{u}_{\varepsilon})\, ds \\
&\leq \|S_x(\varepsilon, \Pi \mathbf{u}_{\varepsilon} )\,\|_{L^2(\Gamma_C)}\|\Pi \mathbf{u}_{\varepsilon, h}- L_h(\Pi \mathbf{u}_{\varepsilon, h})\,\|_{L^2(\Gamma_C)} + \| S_x (\varepsilon, L_h(\Pi \mathbf{u}_{\varepsilon, h}))\,\| _{L^2(\Gamma_C)} \|L_h (\Pi i_h \mathbf{u}_{\varepsilon})- \Pi \mathbf{u}_{\varepsilon}\,\|_{L^2(\Gamma_C)} %\nonumber\\[0.15cm]
+ \alpha \|L_h (\Pi \mathbf{u}_{\varepsilon, h})-\Pi \mathbf{u}_{\varepsilon} \|^2_{L^2(\Gamma_C)}. %\label{eq6}
 \end{align*}
First, by (\ref{ass_3}), we have 
  %\begin{equation}\label{est001}
  \[
  \|S_x(\varepsilon, \Pi \mathbf{u}_{\varepsilon})\|_{L^2(\Gamma_C)}\leq c(1+\|\Pi \mathbf{u}_{\varepsilon}\|_{L^2(\Gamma_C)}) % \quad \forall v\in L^2 (\Gamma_C)
  \]
  %\end{equation}
 and due to (\ref{Gl4}),
 %\begin{equation} \label{est_apriori1_2}
 \[
 \|\Pi \mathbf{u}_{\varepsilon, h}- L_h(\Pi \mathbf{u}_{\varepsilon, h})\,\|_{L^2(\Gamma_C)} \leq c h^{1/2} \|\mathbf{u}_{\varepsilon, h}\|_{\mathbf{H}^{1/2}(\Gamma_C)}.
 %\end{equation}
 \]
 Then, by the triangle inequality, and using the estimates (\ref{Gl4}) and (\ref{est_apriori1}), we have
\begin{align*}%\label{est_apriori2}
\hspace{-0.3cm}\|L_h(\Pi i_h \mathbf{u}_\varepsilon)-\Pi \mathbf{u}_\varepsilon\|_{L^2(\Gamma_C)} \leq \|L_h(\Pi i_h \mathbf{u}_\varepsilon)-\Pi i_h \mathbf{u}_\varepsilon\|_{L^2(\Gamma_C)} + \|\Pi i_h \mathbf{u}_\varepsilon- \Pi \mathbf{u}_\varepsilon\|_{L^2(\Gamma_C)} %\nonumber \\ 
&\leq c_1 h^{1/2} \|i_h \mathbf{u}_\varepsilon\|_{H^{1/2}(\Gamma_C)}+ c_2\|i_h \mathbf{u}_\varepsilon - \mathbf{u}_\varepsilon\|_{H^{1/2}(\Gamma_C)}\\
&\leq  c_1h^{1/2}\left( c h \|\mathbf{u}_\varepsilon\|_{\mathbf{H}^{3/2}(\Gamma_C)} + \|\mathbf{u}_\varepsilon\|_{\mathbf{H}^{1/2}(\Gamma_C)} \right) + \tilde{c}_2h \|\mathbf{u}_\varepsilon\|_{\mathbf{H}^{3/2}(\Gamma_C)}.
\end{align*}
Analogously,
%\begin{equation}\label{est_apriori3}
\[
 \|L_h (\Pi \mathbf{u}_{\varepsilon, h})-\Pi \mathbf{u}_{\varepsilon} \|^2_{L^2(\Gamma_C)} \leq  2\|L_h (\Pi \mathbf{u}_{\varepsilon, h})- \Pi \mathbf{u}_{\varepsilon, h} \|^2_{L^2(\Gamma_C)}+ 2\|\Pi \mathbf{u}_{\varepsilon, h}-\Pi \mathbf{u}_{\varepsilon} \|^2_{L^2(\Gamma_C)} \leq  \tilde{c}_3 h\| \mathbf{u}_{\varepsilon,h}\|^2_{H^{1/2}(\Gamma_C)}+  \tilde{c}_4\|\mathbf{u}_\varepsilon-\mathbf{u}_{\varepsilon,h}\|^2_{H^{1/2}(\Gamma_C)}.
 \]
 % \end{equation}
 %Combining (\ref{est001}) with (\ref{Gl1}) and (\ref{Gl4}), and 
 %using thereby the estimates (\ref{est_apriori1}) - (\ref{est_apriori3}), 
Therefore, we have 
%Hence, using similar arguments as above, see (\ref{f1_1}),  we can conclude
\begin{align*}
D  \leq c_1 h^{1/2}(1+\|\mathbf{u}_\varepsilon\|_{H^{1/2}(\Gamma_C)})\|  \mathbf{u}_{\varepsilon, h}\|_{H^{1/2}(\Gamma_C)} & +c_2 (1+\|\mathbf{u}_{\varepsilon, h}\|_{H^{1/2}(\Gamma_C)})
(c_3h^{3/2}\|\mathbf{u}_\varepsilon\|_{H^{3/2}(\Gamma_C)}+ c_4 h \|\mathbf{u}_\varepsilon\|_{H^{3/2}(\Gamma_C)} + c_5 h^{1/2}\|\mathbf{u}_\varepsilon\|_{H^{1/2}(\Gamma_C)}) \\
%\|\Pi  \mathbf{u}_{\varepsilon}\|_{H^{1/2}(\Gamma_C)}
& +  \alpha \tilde{c}_3   h \| \mathbf{u}_{\varepsilon,h}\|_{H^{1/2}(\Gamma_C)}^2+ \alpha \tilde{c}_4\|\mathbf{u}_\varepsilon-\mathbf{u}_{\varepsilon,h}\|^2_{H^{1/2}(\Gamma_C)}.
\end{align*}
Hence, taking into account the uniform boundedness of $\{\|\mathbf{u}_{\varepsilon,h}\|_{\tilde{\mathbf{H}}^{1/2}(\Gamma_0)}\}$, we conclude that there exists a constant $c=c(\mathbf{u}_\varepsilon)$ such that 
%\begin{equation} \label{Dest}
\[
D \leq c h^{1/2} + \alpha \tilde{c}_2 \|\mathbf{u}_\varepsilon-\mathbf{u}_{\varepsilon,h}\|^2_{\mathbf{H}^{1/2}(\Gamma)}.
\]
%\end{equation}
Altogether yields the claimed estimate of the error $\|\mathbf{u}_{\varepsilon,h}- \mathbf{u}_\varepsilon\|_{\mathbf{H}^{1/2}(\Gamma)}$ 
provided that  
$\alpha $ in (\ref{uniq_reg}) is small enough.
\qed
 \section{Numerical experiments} 
\begin{figure}[t!]
\centering
\includegraphics[trim = 0cm 0.5cm 0cm 0cm, clip,width=1.0\textwidth]{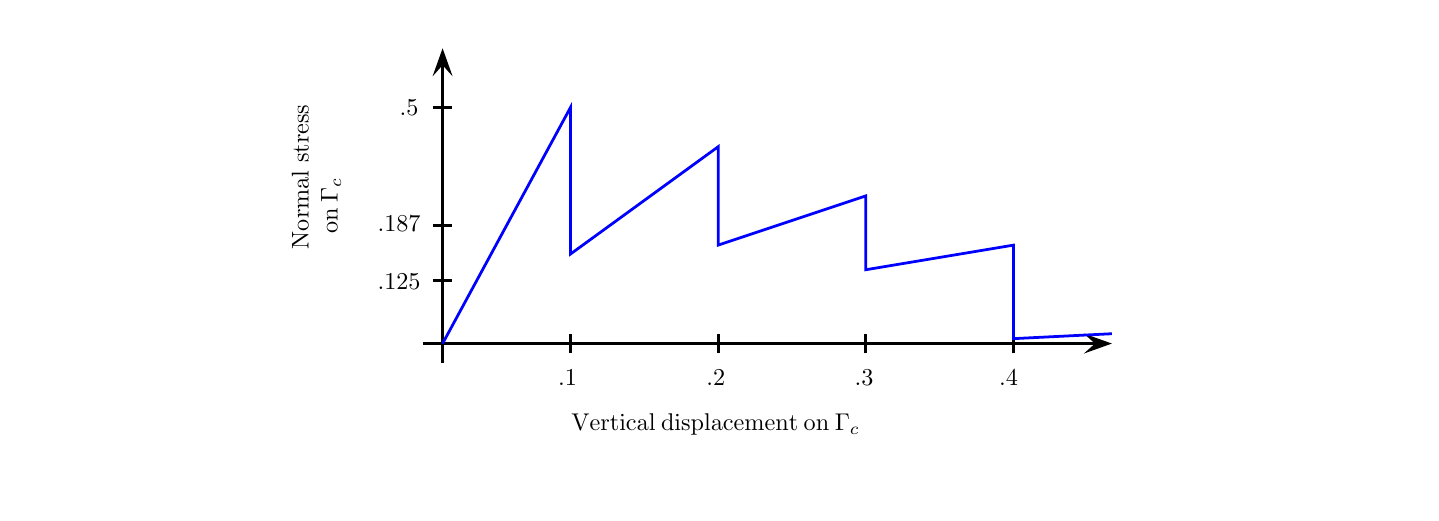}
\caption{Nonmonotone adhesion law $\partial f$ } \label{fig_del}
\end{figure}

For the numerical experiments, we choose $\Omega=(0,100)\times (0,10)$ in [mm]. The boundary is divided into three parts $\Gamma_D=\{0\}\times [0,10]$, $\Gamma_C=(0, 100]\times \{0\}$ and $\Gamma_N=\partial \Omega \backslash (\Gamma_D \cup \Gamma_C)$. The two-dimensional example can
be treated as an approximation for a three-dimensional case considering the domain $\Omega$ as the  cross section of a three-dimensional linear elastic body.  The material parameters are $E=210$ GPa and $\nu= 0.3$. The applied loads on $[50, 100] \times \{10\}$ are $(0,t_2)$ with $t_2=$  $0.2$, $0.4$, $0.6$, $0.8$, $1.0$ $[N/mm^2]$, respectively. The applied loads on $\{100\}\times[100, 10]$ are zero, the volume forces $\mathbf{f}$ are also neglected. We model the nonmonotone adhesion law depicted on Figure \ref{fig_del} with minimum superpotential $f$ defined by 
\begin{eqnarray} \label{min_superpotential}
f(u_n(s))   =  \min\{g_1(-u_n(s)), g_2(-u_n(s)), g_3(-u_n(s)), g_4(-u_n(s)), g_5(-u_n(s))\} \nonumber \\
 = -\max  \{-g_1(-u_n(s)), -g_2(-u_n(s)), -g_3(-u_n(s)), -g_4(-u_n(s)), -g_5(-u_n(s))\}
\end{eqnarray}
with
\[
g_1(y)=\frac{A_1}{2t_1}y^2, \; g_2(y)=b_2(y^2-t_1^2)+d_2, \;g_3(y)=b_3(y^2-t_2^2)+d_3, \; g_4(y)=b_4(y^2-t_3^2)+d_2, \; g_5(y)=d_5
\]
and parameters
\[
A_1=0.5 N/mm^2, \; A_2=0.4375 N/mm^2, \; A_3=0.3125 N/mm^2, \; A_4=0.1875  N/mm^2, 
\]
\[
t_1=0.1 mm, \; t_2=0.2 mm, \; t_3=0.3 mm, \; t_4=0.4 mm, %\mbox{in}\; mm,
\]
\[
b_2=\frac{A_2}{2t_2}, \; d_2= A_1\frac{t_1}{2},\;  b_3=\frac{A_3}{2t_3}, \; d_3=b_2(t_2^2-t_1^2)+d_2, \; b_4=\frac{A_4}{2t_4}, \; d_4=b_3(t_3^2-t_2^2)+d_3, \; d_5=b_4(t_4^2-t_3^2)+d_4.
\]
All computations use piecewise linear functions on a uniform grid with $160$ nodes. The number of the unknowns in the discrete regularized boundary problem (\ref{reg_dis_problem}) is $166$ ($40$ nodes on $\Gamma_C\backslash \overline{\Gamma}_D$ and $43$ nodes on $\Gamma_N$). %All tests have been performed with MATLAB R2013a.  
The regularization parameter $\varepsilon$ is set to $\varepsilon=0.1$. This choice of $\varepsilon$ is based on the observation that smaller values do not improve the solution from a mechanics point of view. The smoothing approximation of the minimum superpotential  (\ref{min_superpotential}) can be handled as above and is based on the smoothing function (\ref{smooth_Zang}). The discrete regularized problem (\ref{reg_dis_problem}) is solved using the following steps. Firstly, we rewrite (\ref{reg_dis_problem}) as a mixed complementarity problem. Secondly, by using the Fischer-Burmeister function $f(a,b)=\sqrt{a^2+b^2}-(a+b)$  we reformulate the last problem  as a system of nonlinear equations of  the form $F(\cdot) = 0$.  Finally, by using $\frac{1}{2}\|F(\cdot)\|^2$ as a merit function,
we obtain a smooth unconstrained minimization problem, which is solved by the {\it lsqnonlin} MATLAB function based on the trust-region Newton method. For details, we
refer the reader to \cite{Ovcharova2012}. The maximal number of iterations in {\it lsqnonlin} has been fixed to $100$. 

\begin{figure}[h!]
\centering
\includegraphics[trim = 0cm 10cm 0cm 9cm, clip,width=1.0\textwidth]{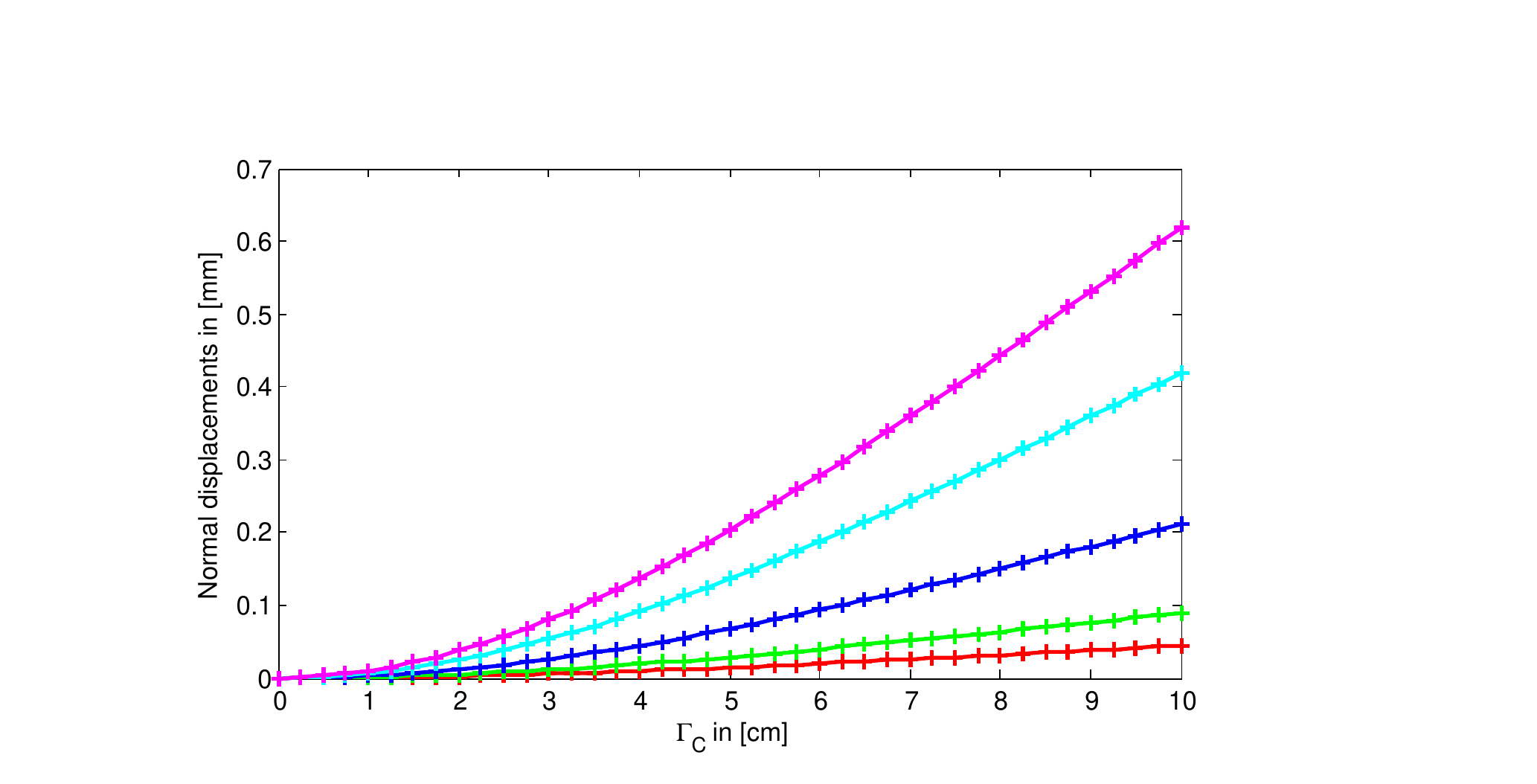}
\caption{The vertical displacements on $\Gamma_C$ for 5 different loads } \label{fig_displ}
\end{figure}

\begin{figure}[h!]
\centering
\includegraphics[trim = 0cm 10cm 0cm 10cm, clip,width=1.0\textwidth]{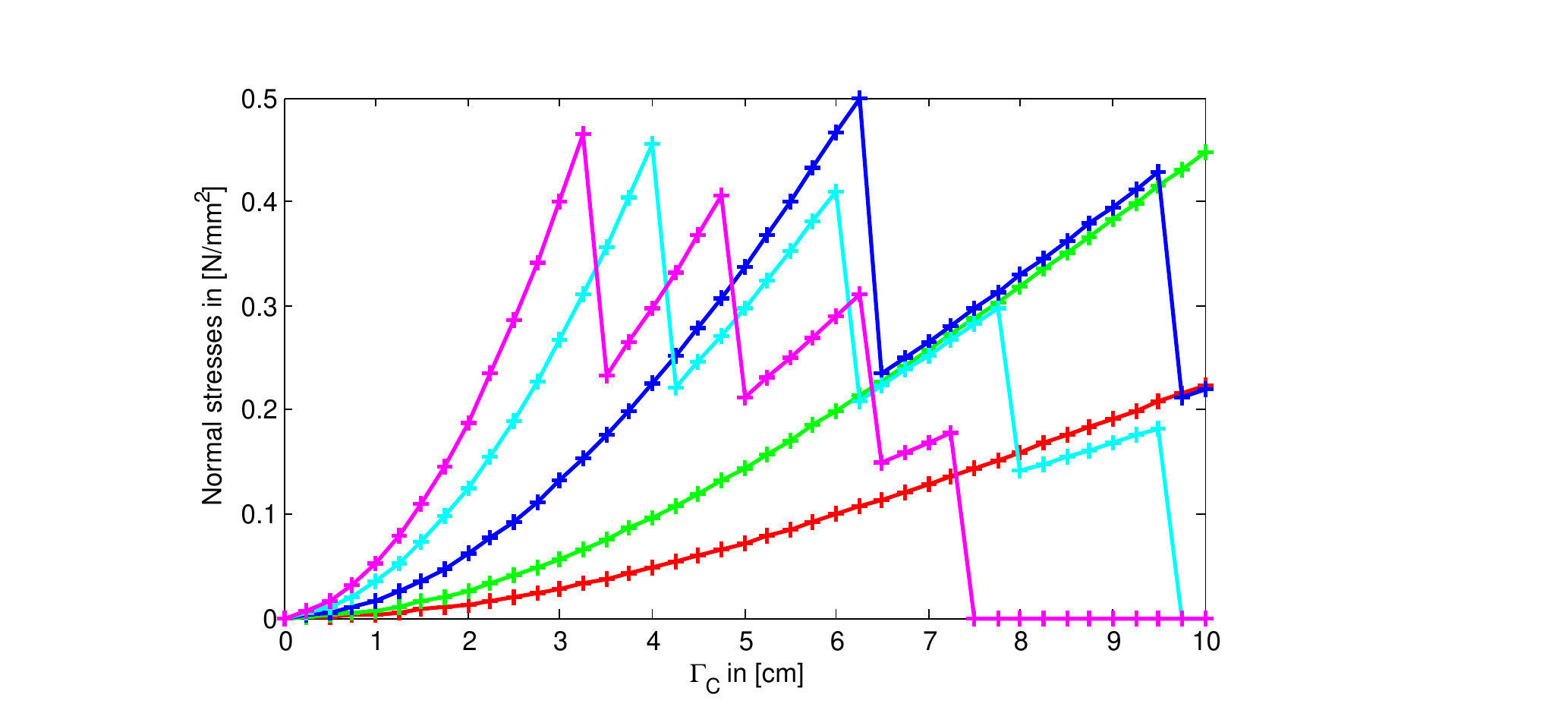}
\caption{The normal stresses on $\Gamma_C$ for 5 different loads} \label{fig_stresses}
\end{figure}

The numerical results are plotted on Figures \ref{fig_displ} and \ref{fig_stresses}. They illustrate the computed vertical displacements %$-\mathbf{u}_n$ 
and the  normal component $\sigma_n$ of the boundary stress vector along $\Gamma_C$. One can see that the computed normal stresses on $\Gamma_C$ reflect the adhesion law from Figure \ref{fig_del}. 
\section* {Conclusions}
%{\color{blue} 
In this paper, we have presented a novel approximation method for solving hemivariational inequalities. This method is based on a smooth approximation of the nonsmooth functional and then, discretization by  $h$-BEM. 
As a future work, we can combine the %smoothing of the data 
the regularization techniques
with  $hp$-adaptive BEM to improve the  convergence
rates of the discretization based on appropriate and automated mesh refinements ($h$-adaptivity) and raising of the polynomial degree ($p$-adaptivity).
%proposed approximation scheme for numerical treatment of boundary hemivariational inequalities.  
Another interesting direction of research is related to the development of the reliable a-posteriori error estimates for the nonsmooth variational problems, which are up to now still missing in the literature.
\\
{\bf Acknowledgment} The author is  grateful to the reviewer for the constructive  comments and suggestions that significantly improved  the manuscript.  

\end{document}